\documentclass[12pt]{amsart}
\usepackage{amssymb}

\newcommand{\x}{\times}
\newcommand{\<}{\langle}
\renewcommand{\>}{\rangle}

\renewcommand{\a}{\alpha}
\renewcommand{\b}{\beta}
\renewcommand{\d}{\delta}
\renewcommand{\o}{\omega}

\newcommand{\e}{\varepsilon}
\newcommand{\g}{\gamma}
\newcommand{\G}{\Gamma}
\renewcommand{\l}{\lambda}

\newcommand{\var}{\varphi}
\newcommand{\s}{\sigma}
\newcommand{\Sig}{\Sigma}
\renewcommand{\t}{\tau}
\renewcommand{\th}{\theta}

\newcommand{\p}{\partial}
\newcommand{\bell}{\mbox{\boldmath$\ell$}}

\newtheorem{theorem}{Theorem}[section]
\newtheorem{mtheorem}{Main Theorem}[section]
\newtheorem{lemma}{Lemma}[section]
\newtheorem{corollary}{Corollary}[section]
\newtheorem{proposition}{Proposition}[section]
\newtheorem{notation}{Notation}[section]
\newtheorem{definition}{Definition}[section]
\newtheorem{question}{Question}[section]

\newtheorem{example}{Example}[section]

\begin{document}

\title[Group $C^*$-Algebras]{Group $C^*$-Algebras as \\ Compact 
Quantum Metric Spaces}
\author{Marc A. Rieffel}
\address{Department of Mathematics \\
University of California \\ Berkeley, CA 94720-3840}
\email{rieffel@math.berkeley.edu}
\date{December 18, 2002}
\thanks{The research reported here was
supported in part by National Science Foundation grant DMS99-70509.}
\subjclass
%[2000]
{Primary 47L87; Secondary 20F65, 53C23, 58B34}

\begin{abstract}
Let $\ell$ be a length function on a group $G$, and let $M_{\ell}$ denote the
operator of pointwise multiplication by $\ell$ on $\bell^2(G)$. 
Following Connes,
$M_{\ell}$ can be used as a ``Dirac'' operator for $C_r^*(G)$.  It defines a
Lipschitz seminorm on $C_r^*(G)$, which defines a metric on the state space of
$C_r^*(G)$.  We investigate whether the topology from this metric 
coincides with the
weak-$*$ topology (our definition of a ``compact quantum metric 
space'').  We give an
affirmative answer for $G = {\mathbb Z}^d$ when $\ell$ is a word-length, or the
restriction to ${\mathbb Z}^d$ of a norm on ${\mathbb R}^d$.  This works for
$C_r^*(G)$ twisted by a $2$-cocycle, and thus for non-commutative tori.  Our
approach involves Connes' cosphere algebra, and an interesting 
compactification of
metric spaces which is closely related to geodesic rays.
\end{abstract}
    
\maketitle
\allowdisplaybreaks

\setcounter{section}{-1}
\section{Introduction}
\label{sec0}

The group $C^*$-algebras of discrete groups provide a much-studied class of
``compact non-commutative spaces'' (that is, unital $C^*$-algebras). 
In \cite{C1}
Connes showed that the ``Dirac'' operator of an unbounded Fredholm 
module over a
unital $C^*$-algebra provides in a natural way a metric on the state 
space of the
algebra.  Unbounded Fredholm modules (i.e. spectral triples)
also provide smooth structure, important
homological data and much else.  In the subsequent years Connes has 
been strongly
advocating this use of Dirac operators as the way to deal with the Riemannian
geometry of non-commutative spaces \cite{C2}, \cite{Cn3}, \cite{C5}, 
\cite{C4}.  The
class of examples most discussed in \cite{C1} consists of the  group 
$C^*$-algebras
of discrete groups, with the Dirac operator coming in a simple way 
from a length
function on the group.  Connes obtained in \cite{C1} strong 
relationships between
the growth of a group and the summability of Fredholm modules over its group
$C^*$-algebra.  However he did not explore much the metric on the state space.

In \cite{R4}, \cite{R5} I pointed out that, motivated by what happens 
for ordinary
compact metric spaces, it is natural to desire that the topology from 
the metric on
the state space coincides with the weak-$*$ topology (for which the 
state space is
compact).  This property was verified in \cite{R4} for certain 
examples, notably the
non-commutative tori, with ``metric'' structure coming from a different
construction.   (See \cite{R5}, \cite{R6}, \cite{R7} for further 
developments.)  But
in general I have found this property to be difficult to verify for 
many natural
examples.

The main purpose of this paper is to examine this property for 
Connes' initial class
of examples, the group $C^*$-algebras with the Dirac operator coming 
from a length
function.  To be more specific, let $G$ be a countable (discrete) 
group, and let
$C_c(G)$ denote the convolution $*$-algebra of complex-valued 
functions of finite
support on $G$.  Let $\pi$ denote the usual $*$-representation of $C_c(G)$ on
$\bell^2(G)$ coming from the unitary representation of $G$ by left 
translation on
$\bell^2(G)$.  The norm-completion of $\pi(C_c(G))$ is by definition 
the reduced
group $C^*$-algebra, $C_r^*(G)$, of $G$.  We identify $C_c(G)$ with 
its image in
$C_r^*(G)$, so that it is a dense $*$-subalgebra.

Let a length function $\ell$ be given on $G$.  We let $M_{\ell}$ 
denote the (usually
unbounded) operator on $\bell^2(G)$ of pointwise multiplication by 
$\ell$.  Then
$M_{\ell}$ will serve as our ``Dirac'' operator.  One sees easily 
\cite{C1} that the
commutators $[M_{\ell},\pi_f]$ are bounded operators for each $f \in 
C_c(G)$.  We
can thus define a seminorm, $L_{\ell}$, on $C_c(G)$ by $L_{\ell}(f) =
\|[M_{\ell},\pi_f]\|$.

In general, if $L$ is a seminorm on a dense $*$-subalgebra $A$ of a unital
$C^*$-algebra ${\bar A}$ such that $L(1) = 0$, we can define a 
metric, $\rho_L$, on
the state space $S({\bar A})$ of ${\bar A}$, much as Connes did, by
\[
\rho_L(\mu,\nu) = \sup\{|\mu(a) - \nu(a)|: a \in A,\ L(a) \le 1\}.
\]
(Without further hypotheses $\rho_L$ may take value $+\infty$.)  In 
\cite{R5} we
define
$L$ to be a {\em Lip-norm} if the topology on $S({\bar A})$ from 
$\rho_L$ coincides
with the weak-$*$ topology.  We consider a unital $C^*$-algebra 
equipped with a
Lip-norm to be a compact quantum metric space \cite{R5}.

The main question dealt with in this paper is whether the seminorms $L_{\ell}$
coming as above from length functions on a group are Lip-norms.  In 
the end we only
have success in answering this question for the groups ${\mathbb Z}^d$.  The
situation there is already somewhat complicated because of 
the large variety of
possible length-functions.  But we carry out our whole discussion in 
the slightly
more general setting of group $C^*$-algebras twisted by a $2$-cocycle 
(definitions
given later), and so this permits us to treat successfully also the 
non-commutative
tori \cite{R3}.  The main theorem of this paper is:

\begin{mtheorem}
\label{th0.1}
Let $\ell$ be a length function on ${\mathbb Z}^d$ which is either 
the word-length
function for some finite generating subset of ${\mathbb Z}^d$, or the 
restriction to
${\mathbb Z}^d$ of some norm on ${\mathbb R}^d$.  Let $c$ be a $2$-cocycle on
${\mathbb Z}^d$, and let $\pi$ be the regular representation of $C^*({\mathbb
Z}^d,c)$ on $\bell^2({\mathbb Z}^d)$.  Then the seminorm $L_{\ell}$ defined on
$C_c({\mathbb Z}^d)$ by $L_{\ell}(f) = \|[M_{\ell},\pi_f]\|$ is a Lip-norm on
$C^*({\mathbb Z}^d,c)$.
\end{mtheorem}

The path which I have found for the proof of this theorem is somewhat 
long, but it
involves some objects which are of considerable independent interest, 
and which may
well be useful in treating more general groups.  Specifically, we 
need to examine
Connes' non-commutative cosphere algebra \cite{C5} for the examples which we
consider.  This leads naturally to a certain compactification which one can
construct for any locally compact metric space.  We call this ``the metric
compactification''.  Actually, this compactification had been introduced 
much earlier by
Gromov \cite{G3}, but it is different from the famous Gromov 
compactification for a
hyperbolic metric space, and it seems not to have received much 
study.  Our approach
gives a new way of defining this compactification.  
We also need to examine the
strong relationship between geodesic rays and points in the boundary of this
compactification, since this will provide us with enough points of 
the boundary which have finite orbits.  For word-length functions 
on ${\mathbb Z}^d$ this is already
fairly complicated.

The contents of the sections of this paper are as follows.  In 
Section~1 we make
more precise our notation, and we make some elementary observations 
showing that on
any separable unital $C^*$-algebra there is an abundance of Lip-norms, 
and that
certain constructions in the literature concerning groups of ``rapid 
decay''  yield
natural Lip-norms on $C^*_r(G)$.  In Section~2 we begin our 
investigation of the
Dirac operators for $C_r^*(G)$ coming from length functions.  In 
Section~3 we examine
Connes' cosphere algebra for our situation. We 
show in particular that if
the action of the group on the boundary of its metric compactification is 
amenable, then
the cosphere algebra has an especially simple description.  Then in 
Section~4 we study the metric compactification in general,
with attention to the geodesic rays.

In Section~5 we begin our study of specific groups by considering the group
${\mathbb Z}$.  This is already interesting.  (Consider a generating 
set such as
$\{\pm 3,\pm 8\}$.)  The phenomena seen there for ${\mathbb Z}$
indicate some of the complications
which we will encounter in trying to deal with ${\mathbb Z}^d$.  In 
Section~6 we
study the metric compactification of ${\mathbb R}^d$ for any given 
norm, and then in
Section~7 we apply this to prove the part of our Main Theorem for 
length functions
on ${\mathbb Z}^d$ which are the restrictions of norms on ${\mathbb R}^d$.  In
Section~8 we study the metric compactification of ${\mathbb Z}^d$ for 
word-length
functions, and in Section~9 we apply this to prove the remaining part 
of our Main
Theorem.  We conclude in Section~10 with a brief examination of the free
(non-Abelian) group on two generators, to see both how far our 
approach works, and
where we become blocked from proving for it the corresponding version 
of our Main Theorem.

Last-minute note: I and colleagues believe we have a proof that the
Main Theorem is also true for word-hyperbolic groups with word-length functions,
using techniques which are entirely different from those used here,
and which do not seem to apply to the case of ${\mathbb Z}^d$ treated
here.

A substantial part of the research reported here was carried out while I
visited the Institut de Math\'ematique de Luminy, Marseille, for three
months.  I would like to thank Gennady Kasparov, Etienne Blanchard, Antony
Wasserman, and Patrick Delorme very much for their warm hospitality and
their mathematical stimulation during my very enjoyable visit.

%%%%%%%%%%%%%%%%%%%%%%%%%%%%%%%%%%%%%%%%%%%%%%%%
\setcounter{section}{0}
\section{An abundance of Lip-norms}
\label{sec1}
%%%%%%%%%%%%%%%%%%%%%%%%%%%%%%%%%%%%%%%%%%%%%%%%

In this section we establish some of our notation, and show that on 
any separable
unital $C^*$-algebra there is an abundance of Lip-norms.  In the 
absence of further
structure these Lip-norms appear somewhat artificial.  But we then 
show that some
known constructions for group $C^*$-algebras yield somewhat related 
but more natural
Lip-norms.

Our discussion in the next few paragraphs works in the greater generality of
order-unit spaces which was used in \cite{R5}.  But we will not use 
that generality
in later sections, and so the reader can have in mind just the case 
of dense unital *-subalgebras of unital
$C^*$-algebras, with the identity element being the order unit.  We 
recall that a (possibly discontinuous)
seminorm $L$ on an order-unit space is said to be lower 
semicontinuous if $\{a \in A:
L(a) \le r\}$ is norm-closed for any $r > 0$.
        
\begin{proposition}
\label{prop1.1}
Let $A$ be an order-unit space which is separable.  For any countable 
subset $E$ of
$A$ there are many lower semicontinuous Lip-norms on $A$ which are defined and
finite on $E$.
\end{proposition}

\begin{proof}
The proof is a minor variation on the fact that the weak-$*$ topology 
on the unit
ball of the dual of a separable Banach space is metrizable (theorem V.5.1 of
\cite{Cw}).  We scale each non-zero element of $E$ so that it is in 
the unit ball of
$A$ (and $\ne 0$), and we incorporate $E$ into a sequence, $\{b_n\}$, 
of elements of
$A$ which is dense in the unit ball of $A$.  Let $\{\o_n\}$ be any sequence in
${\mathbb R}$ such that $\o_n > 0$ for each $n$ and $\Sig \o_n < 
\infty$.  Define a
norm, $M$, on the dual space $A'$ of $A$ by
\[
M(\l) = \Sig \o_n |\l(b_n)|.
\]
The metric from this norm, when restricted to the unit ball of $A'$, gives the
weak-$*$ topology, because it is easily checked that if a net in the 
unit ball of
$A'$ converges for the weak-$*$ topology then it converges for the 
metric from $M$,
and then we can apply the fact that the unit ball is weak-$*$ compact.

We let $S(A)$ denote the state space of $A$.  Since $S(A)$ is a 
subset of the unit ball of $A'$, the restriction to $S(A)$ of the 
metric from the
norm $M$ gives $S(A)$ the weak-$*$ topology.  Let $L_M$ denote the 
corresponding
Lipschitz seminorm on $C(S(A))$ from this metric, 
allowing value $+\infty$.  View each 
element of $A$
as a function on $S(A)$ in the usual way.  Then $L_M(b_n) \le 
\o_n^{-1} < \infty$ for each
$n$, because if $\mu,\nu \in S(A)$ then
\[
|b_n(\mu) - b_n(\nu)| = |(\mu-\nu)(b_n)| \le \o_n^{-1}M(\mu-\nu) = \o_n^{-1}
\rho_M(\mu,\nu).
\]
Let $B$ denote the linear span of $\{b_n\}$ together with the order 
unit.  Then $B$
is a dense subspace of $A$ containing the order-unit, and $L_M$ 
restricted to $B$ is
a seminorm which can be verified to be lower semicontinuous.  The 
inclusion of $A$
into $C(S(A))$ is isometric (on self-adjoint elements if $A$ is a 
$C^*$-algebra) and
since $L_M$ comes from an ordinary metric, it follows that $L_M$ on $A$ is a
Lip-norm.  (For example, use theorem $1.9$ of \cite{R4}.)
\end{proof}

The considerations above are close to those of theorem $9.8$ of 
\cite{R5}.  Let me
take advantage of this to mention here that Hanfeng Li showed me by clever
counterexample that theorem $9.8$ of \cite{R5} is not correct as 
presented, because
$A$ may not be big enough.  However, if $A$ is taken to be 
norm-complete, then there
is no difficulty.  Theorem $9.11$ needs to be adjusted accordingly. 
But this change
does not affect later sections of \cite{R5} nor the subsequent papers 
\cite{R6},
\cite{R7}.

We now turn to (twisted) group $C^*$-algebras, 
and we use a different approach,
which takes advantage of the fact that the group elements provide a 
natural ``basis''
for the group $C^*$-algebras.  Thus let $G$ be a countable discrete 
group, and let
$c$ be a $2$-cocycle \cite{ZM} on $G$ with values 
in the circle group ${\mathbb
T}$.  We assume that $c$ is normalized so that $c(x,y) = 1$ if $x = 
e$ or $y = e$.
We let $C^*(G,c)$ denote the full $c$-twisted group $C^*$-algebra of 
$G$, and we
let $C_r^*(G,c)$ denote the reduced $c$-twisted group $C^*$-algebra \cite{ZM},
\cite{Pdr} coming from the left regular representation, $\pi$, on 
$\bell^2(G)$.    Both $C^*$-algebras are completions 
of $C_c(G)$,
the space of finitely supported ${\mathbb C}$-valued functions on $G$, with
convolution twisted by $c$.  Our conventions, following \cite{ZM}, are that
\begin{eqnarray*}
(f*g)(x) &= &\Sig f(y)g(y^{-1}x)c(y,y^{-1}x), \\
f^*(x) &= &{\bar f}(x^{-1}){\bar c}(x,x^{-1}).
\end{eqnarray*}
The left regular representation is given by 
the same formula as the
above twisted convolution, but with $g$ viewed as an element
of $\bell^2(G)$. Then $C_r^*(G,c)$ is the completion of $C_c(G)$
for the operator norm coming from the left regular representation. 
We will often set $A = C_r^*(G,c)$. We note that when $\pi$ is
restricted to $G$ we have
\[
(\pi_y\xi)(x) = \xi(y^{-1}x)c(y,y^{-1}x)
\]
for $\xi \in \bell^2(G)$ and $x,y \in G$.  In particular $\pi_y\pi_z
= c(y,z)\pi_{yz}$.

There is a variety of norms on $C_c(G)$ which have been found to be useful in
addition to the $C^*$-norms.  These other norms are not necessarily 
algebra norms.
To begin with, there is the $\bell^1$-norm, as well as the 
$\bell^p$-norms for $1 < p
\le \infty$.  But let $\ell$ be a length function on $G$, so that 
$\ell(xy) \le \ell(x)
+ \ell(y)$, $\ell(x^{-1}) = \ell(x)$, $\ell(x) \ge 0$, and $\ell(x) = 
0$ exactly if
$x = e$, the identity element of $G$.  Then in connection with groups 
of ``rapid
decay'' (such as word-hyperbolic groups) one defines norms on $C_c(G)$ of the
following form \cite{Jo1}, \cite{Jo2}, \cite{Ji}, \cite{JiS}:
\[
\|f\|_{p,k} = (\Sig (|f(x)|(1+\ell(x))^k)^p)^{1/p}.
\]
These norms clearly have the properties that
\begin{itemize}
\item[1)] $\|f\|_{p,k} \le \||f|\|_{p,k}$ (actually $=$),
\item[2)] if $|f| \le |g|$ then $\||f|\|_{p,k} \le \||g|\|_{p,k}$.
\end{itemize}
Their interest lies in the fact that for 
a rapid-decay group and an appropriate
choice of $p$ and $k$ depending on the group, one has (see the first 
line of the
proof of theorem $1.3$ of \cite{Ji}, combined, in the case of 
nontrivial cocycle,
with proposition $3.10$b of \cite{JiS}):
\begin{itemize}
\item[3)] There is a constant, $K$, such 
that $\|f\|_{C_r^*} \le K\|f\|_{p,k}$.
\end{itemize}
Notice also that if the cocycle $c$ is trivial and if $G$ is amenable 
\cite{Pdr} then
the $C^*$-norm itself satisfies the above three properties, because 
from the trivial
representation we see that for $f \in C_c(G)$ we have
\[
\|f\|_{C^*(G)} \le \|f\|_1 = \||f|\|_{C^*(G)},
\]
while if $|f| \le |g|$ then
\[
\||f|\|_{C^*(G)} = \|f\|_1 \le \|g\|_1 = \||g|\|_{C^*(G)}.
\]
Finally, for any group and any cocycle we always have at least the 
$\bell^1$-norm
which satisfies the above three properties.

With these examples in mind, we make

\setcounter{definition}{1}
\begin{definition}
\label{def1.2}
Let $\|\cdot\|_A$ denote the $C^*$-norm on $A = C_r^*(G,c)$.  We will 
say that a
norm, $\|\cdot\|$, on $C_c(G)$ is {\em order-compatible} with 
$\|\cdot\|_A$ if for
all $f,g \in C_c(G)$ we have:
\begin{itemize}
\item[1)] $\|f\| \le \||f|\|$.
\item[2)] If $|f| \le |g|$ then $\||f|\| \le \||g|\|$.
\item[3)] There is a constant, $K$, such that $\|f\|_A \le K\|f\|$.
\end{itemize}
\end{definition}

We remark that these conditions are a bit weaker than those
required for a ``good norm'' in \cite{Lf}.

Suppose now that $\o$ is a real-valued function on $G$ such that 
$\o(e) = 0$ and
$\o(x) > 0$ for $x \ne e$.  Fix an order-compatible norm $\|\cdot\|$ 
on $C_c(G)$,
and set
\[
L(f) = \|\o|f|\|.
\]
It is clear that $L$ is a seminorm which is $0$ only on the span of 
the identity
element of the convolution algebra $C_c(G,c)$.  (Thus $L$ is a 
Lipschitz seminorm as
defined in \cite{R5}.)  In the way discussed in the introduction, 
$L$ defines a
metric, $\rho_L$, on $S(C_r^*(G,c))$ by
\[
\rho_L(\mu,\nu) = \sup\{|\mu(f) - \nu(f)|: L(f) \le 1\},
\]
which may take value $+\infty$.  Denote $C_r^*(G,c)$ by $A$, and its
C*-norm by
$\|\cdot\|_A$, as above.

\setcounter{lemma}{2}
\begin{lemma}
\label{lem1.3}
Suppose that there is a constant $s > 0$ such that $\o(x) \ge s$ for 
all $x \ne e$.
Then $\rho_L$ gives $S(A)$ finite radius.  (In particular, $\rho_L$ 
does not take
the value $+\infty$.)
\end{lemma}

\begin{proof}
Let $f \in C_c(G)$, and assume that $f(e) = 0$.  Let $K$ be the constant in the
definition of ``order-compatible''.  Then
\[
\|f\|_A \le K\|f\| \le K\||f|\| \le Ks^{-1}\|\o|f|\| = Ks^{-1}L(f).
\]
The desired conclusion then follows from proposition $2.2$ of \cite{R5}.
\end{proof}

\setcounter{lemma}{3}
\begin{lemma}
\label{lem1.4}
Suppose that $\o(x) = 0$ only if $x = e$ and that the function $\o$ 
is ``proper'',
in the sense that for any $n$ the set $\{x \in G: \o(x) \le n\}$ is 
finite (so, in
particular, there exists a constant $s$ as in the above lemma).  Then 
the topology
from the metric $\rho_L$ on $S(A)$ coincides with the weak-$*$ 
topology.  Thus $L$
is a Lip-norm.
\end{lemma}

\begin{proof}
We apply theorem $1.9$ of \cite{R4}.  As in that theorem, we set
\[
{\mathcal B}_1 = \{f \in C_c(G): \|f\|_A \le 1 \mbox{ and } L(f) \le 1\}.
\]
The theorem tells us that it suffices to show that ${\mathcal B}_1$ is totally
bounded for $\|\cdot\|_A$.  So let $\e > 0$ be given.  Adjust $K$ if 
necessary so
that $K \ge 1$, and set
\[
E = \{x \in G: \o(x) \le 3K/\e\}.
\]
Then $E$ is a finite set because $\o$ is proper.  Set $A^E = \{f \in 
C_c(G): f(x) =
0 \mbox{ for } x \notin E\}$, so that $A^E$ is a finite-dimensional subspace of
$C_c(G)$.  In particular, $A^E \cap {\mathcal B}_1$ is totally bounded.

Let $f \in {\mathcal B}_1$.  Then $f = g+h$ where $g \in A^E$ and 
$h(x) = 0$ for $x
\in E$.  Now $|h| \le |f|$, and $\o(x) \ge 3K/\e$ on the support of $h$, and so
\begin{eqnarray*}
\|h\|_A &\le &K\|h\| \le K\||h|\| \le K(\e/3K)\|\o|h|\| \\
&\le &(\e/3)\|\o|f|\| = (\e/3)L(f) \le \e/3.
\end{eqnarray*}
Thus $\|f-g\|_A = \|h\|_A \le \e/3$.  In particular, $\|g\|_A \le 1 + 
(\e/3)$.  Note
also that $L(g) = \|\o|g|\| \le \|\o|f|\| = L(f) \le 1$.  Thus upon 
scaling $g$ by
$(1+\e/3)^{-1}$ if necessary to obtain an element of ${\mathcal 
B}_1$, we see that
$f$ is within distance $2\e/3$ of ${\mathcal B}_1 \cap A^E$.  Thus a 
finite subset
of ${\mathcal B}_1 \cap A^E$ which is $\e/3$ dense in ${\mathcal B}_1 
\cap A^E$ will
be $\e$-dense in ${\mathcal B}_1$.
\end{proof}

\setcounter{lemma}{4}
\begin{lemma}
\label{lem1.5}
Even without $\o$ being proper, or satisfying the condition of Lemma 
\ref{lem1.3},
the seminorm $L$ is lower semicontinuous (with respect to $\|\cdot\|_A$).
\end{lemma}

\begin{proof}
Let $\{f_n\}$ be a sequence in $C_c(G)$ which converges to $g \in C_c(G)$ for
$\|\cdot\|_A$, and suppose that there is an $r \in {\mathbb R}$ such 
that $L(f_n)
\le r$ for all $n$.  Now $\pi_f\d_0 = f$ where on the 
right $f$ is viewed as an
element of $\bell^2$ and $\d_0$ is the ``delta-function'' at $0$.  
Consequently
$\|f\|_A \ge \|f\|_2 \ge \|f\|_{\infty}$.  Thus $f_n$ converges 
uniformly on $G$ to
$g$.  Let $S$ denote the support of $g$, and 
let $\chi_S$ be its characteristic
function.  Then the sequence $\o\chi_S|f_n|$ converges uniformly to 
$\o|g|$.  But
all norms on a finite-dimensional vector space are equivalent, and so
$\o\chi_S|f_n|$ converges to $\o|g|$ for $\|\cdot\|$.  This says that 
$L(\chi_Sf_n)$
converges to $L(g)$.  But $L(\chi_Sf_n) = \|\o\chi_Sf_n\| \le L(f) 
\le r$.  Thus
$L(g) \le r$.
\end{proof}

We combine the above lemmas to obtain:

\setcounter{proposition}{5}
\begin{proposition}
\label{prop1.6}
Let $\o$ be a proper non-negative function on $G$ such that $\o(x) = 
0$ exactly if
$x = e$.  Let $\|\cdot\|$ be an order-compatible norm on $C_c(G)$, and set
\[
L(f) = \|\o|f|\|
\]
for $f \in C_c(G)$.  Then $L$ is a lower 
semicontinuous Lip-norm on $C_r^*(G)$.
\end{proposition}

We remark that when $\o$ is a length function on $G$ and when $\|\cdot\| =
\|\cdot\|_1$, it is well-known and easily seen that $L$ satisfies the 
Leibniz rule
with respect to $\|\cdot\|_1$, that is
\[
L(f*g) \le L(f)\|g\|_1 + \|f\|_1L(g).
\]
But there seems to be no reason why many of the above Lip-norms 
should satisfy the
Leibniz rule with respect to $\|\cdot\|_A$.  And it is not clear to me what
significance the Leibniz rule has for the metric properties which we 
are examining.

%%%%%%%%%%%%%%%%%%%%%%%%%%%%%%%%%%%%%%%%%%%%%%%%%%%%%%
\setcounter{section}{1}
\section{Dirac operators from length functions}
\label{sec2}
%%%%%%%%%%%%%%%%%%%%%%%%%%%%%%%%%%%%%%%%%%%%%%%%%%%%%%%%

In this section we make various preliminary observations about the 
seminorms $L$ which come from using length functions on a group
as ``Dirac'' operators, as described in the 
introduction.  We
also reformulate our main question as concrete questions concerning $C_r^*(G)$
itself.

We use the notation of the previous section, and we let $M_{\ell}$ denote the
(usually unbounded) operator on $\bell^2(G)$ of 
pointwise multiplication by the
length function $\ell$.  We recall from \cite{C1} why the commutators
$[M_{\ell},\pi_f]$ are bounded for $f \in C_c(G)$.  Let $y \in G$ and $\xi \in
\bell^2(G)$.  Then we quickly calculate that
\[
([M_{\ell},\pi_y]\xi)(x) = (\ell(x) - \ell(y^{-1}x))\xi(y^{-1}x)c(y,y^{-1}x).
\]
 From the triangle inequality for $\ell$ we know that $|\ell(x) - 
\ell(y^{-1}x)| \le
\ell(y)$, and so $\|[M_{\ell},\pi_y]\| \le \ell(y)$.  
In fact, this observation
indicates the basic property of $\ell$ which we need for the 
elementary part of our
discussion, namely that, although $\ell$ is usually unbounded, it 
differs from any
of its left translates by a bounded function.

This suggests that we work in the more general context of functions 
having just this
latter property, as this may clarify some aspects.  Additional 
motivation for doing
this comes from the importance which Connes has demonstrated for examining the
effect of automorphisms of the $C^*$-algebra as gauge transformations, and the
resulting effect on the metric.  In Connes' approach the inner 
automorphisms play a
distinguished role, giving ``internal fluctuations'' of the metric \cite{CC1},
\cite{CC2} (called ``internal perturbations'' in \cite{Cn3}).  However, in our
setting we usually do not have available 
the ``first order'' condition which is
crucial in Connes' setting.  We discuss this briefly at the end of 
this section.

Anyway, in our setting the algebra $C_r^*(G,c)$ has some special inner
automorphisms, namely those coming from the elements of $G$.  The automorphism
corresponding to $z \in G$ is implemented on $\bell^2(G)$ by 
conjugating by $\pi_z$. 
When this
automorphism is composed with the representation, the effect is to change $D =
M_{\ell}$ to $M_{\a_z(\ell)}$, where $\a_z(\ell)$ 
denotes the left translate of
$\ell$ by $z$.  But $\a_z(\ell)$ need not again be a length function, 
although it is
translation bounded.  (In order to try to clarify contexts, 
we will from now on
systematically use $\a$ to denote ordinary left translation of 
functions, especially
when those functions are not to be viewed as being in $\bell^2(G)$. 
Our convention
is that $(\a_z\ell)(x) = \ell(z^{-1}x)$.)  We will make frequent use  of the
easily-verified commutation relation that
\[
\pi_yM_h = M_{\a_y(h)}\pi_y
\]
for any function $h$ on $G$ and any $y \in G$, as long as the  domains of 
definitions of
the product operators are respected.  
This commutation relation is what we used
above to obtain the stated fact about the effect of inner automorphisms.

In what follows we will only use real-valued functions to define our Dirac
operators, so that the latter are self-adjoint.  But much of what follows
generalizes easily to complex-valued functions, or to functions with values in
$C^*$-algebras such as Clifford algebras.  These generalizations 
deserve exploration.

To formalize our discussion above we make:

\setcounter{definition}{0}
\begin{definition}
\label{def2.1}
We will say that a (possibly unbounded) real-valued function, $\o$, 
on $G$ is (left)
{\em translation-bounded} if $\o - \a_y\o$ is a bounded function for 
every $y \in
G$.  For $y \in G$ we set $\var_y = \o - \a_y(\o)$.  So the context 
must make clear
what $\o$ is used to define $\var$.  For each $y \in G$ we set $\ell^{\o}(y) =
\|\var_y\|_{\infty}$.
\end{definition}

Thus every length-function on $G$ is translation-bounded.  Any group 
homomorphism
from $G$ into ${\mathbb R}$ is translation bounded.  (E.g., the 
homomorphism $\o(n)
= n$ from ${\mathbb Z}$ to ${\mathbb R}$ which is basically the 
Fourier transform of
the usual Dirac operator on ${\mathbb T}$.)  Linear combinations of
translation-bounded functions are translation bounded.  In particular, 
the sum of a
translation-bounded function with any bounded function is translation 
bounded.  (As
a more general context one could consider any faithful unitary representation
$(\pi,{\mathcal H})$ of $G$ together with an unbounded self-adjoint 
operator $D$ on
${\mathcal H}$ such that $D - \pi_zD\pi_z^*$ is densely defined and 
bounded for each
$z \in G$, and $D$ satisfies suitable non-triviality conditions.  Our later
discussion will indicate why one may also want to require that the
$(\pi_zD\pi_z^*)$'s all commute with each other.)

It is simple to check that the $\var_y$'s satisfy the $1$-cocycle identity
\setcounter{equation}{1}
\begin{equation}
\label{eq2.2}
\var_{yz} = \var_y + \a_y(\var_z).
\end{equation}
We will make use of this relation a number of times.  This type of 
relation occurs
in various places in the literature in connection with dynamical systems.

Simple calculations show that $\ell^{\o}$ satisfies the axioms for a length
function except that we may have $\ell^{\o}(x) = 0$ for some $x \ne e$.  Notice
also that if $\o$ is already a length function, then $\ell^{\o} = \o$.  We also
remark that in general we can always add a constant function to $\o$ without
changing the corresponding $\var_y$'s, $\ell^{\o}$, or the commutators
$[M_{\ell},\pi_y]$.  In particular, we can always adjust $\o$ in this 
way so that
$\o(e) = 0$ if desired.

We now fix a translation-bounded function, $\o$, on $G$, and we consider the
operator, $M_{\o}$, of pointwise multiplication on $\bell^2(G)$.  It is
self-adjoint.  We use it as a ``Dirac operator''.  
The calculation done earlier
becomes
\[
[M_{\o},\pi_y] = M_{\var_y}\pi_y.
\]
 From this we see that for each $y \in G$ we have
\[
\|[M_{\o},\pi_y]\| = \ell^{\o}(y).
\]
For any $f \in C_c(G)$ we have
\[
[M_{\o},\pi_f] = \Sig f(y)M_{\var_y}\pi_y,
\]
and consequently we have
\[
\|[M_{\o},\pi_f]\| \le \|\ell^{\o}f\|_1,
\]
where $\ell^{\o}f$ denotes the pointwise product.  We set
\[
L^{\o}(f) = \|[M_{\o},\pi_f]\|.
\]
Then $L^{\o}$ is a seminorm on $C_c(G) \subseteq C_r^*(G,c)$, and 
$L^{\o}$ is lower
semicontinuous by proposition $3.7$ of \cite{R5}.  A calculation 
above tells us that
$L^{\o}(\d_x) = \ell^{\o}(x)$ for all $x \in G$.  In particular, 
$L^{\o}(\d_e) =
0$, with $\d_e$ the identity element of the convolution algebra $C_c(G)$.

If we view $\d_z$ as the usual basis element at $z$ for $\bell^2(G)$, 
then for any
$f \in C_c(G)$ we have
\[
[M_{\o},\pi_f]\d_z = \Sig f(y)M_{\var_y}c(y,z)\d_{yz}
\]
for each $z$.  From this we easily obtain:

\setcounter{proposition}{2}
\begin{proposition}
\label{prop2.3}
Let $f \in C_c(G)$.  Then $L^{\o}(f) = 0$ exactly if $\var_y = 0$ for 
each $y$ in
the support of $f$, that is, exactly if $\ell^{\o}f = 0$.  Thus if 
$\ell^{\o}(x) >
0$ for all $x \ne e$, then $L^{\o}$ is a Lipschitz seminorm in the 
sense that its
null space is spanned by $\d_e$.
\end{proposition}

We would like to know when $L^{\o}$ is a Lip-norm.  Of course, 
$L^{\o}$ defines, as
earlier, a metric on the state space $S(C_r^*(G,c))$, which may take value
$+\infty$.  We denote this metric by $\rho_{\o}$.  As a first step, 
we would like to
know whether $\rho_{\o}$ gives $S(C_r^*(G,c))$ finite radius.  We recall from
proposition $2.2$ of \cite{R5} that 
this will be the case if there is an $r \in
{\mathbb R}$ such that $\|f\|^{\sim} \le rL(f)$ for all $f \in C_c(G)$, where
$\|f\|^{\sim} = \inf\{\|f - \a\d_e\|: \a \in {\mathbb C}\}$. 
Officially speaking
we should work with self-adjoint $f$'s, but by the comments before 
definition $2.1$
of \cite{R6} we do not need to make this restriction because clearly 
$L^{\o}(f^*) =
L^{\o}(f)$ for each $f$.  However we find it convenient to use the following
alternative criterion for finite radius, which is natural in our 
situation because
we have a canonical tracial state:

\setcounter{proposition}{3}
\begin{proposition}
\label{prop2.4}
Let $L$ be a Lipschitz seminorm on an order-unit space $A$, and let 
$\mu$ be a state
of $A$.  If the metric $\rho_L$ from $L$ gives $S(A)$ finite radius 
$r$, then $\|a\|
\le 2rL(a)$ for all $a \in A$ such that $\mu(a) = 0$.  Conversely, if 
there is a
constant $k$ such that
\[
\|a\| \le kL(a)
\]
for all $a \in A$ such that $\mu(a) = 0$, then $\rho_L$ gives $S(A)$ radius no
greater than $k$.
\end{proposition}

\begin{proof}
Suppose the latter condition holds.  For any given $a \in A$ set $b = a - 
\mu(a)e$.  (Here
$e$ is the order-unit.)  Then $\mu(b) = 0$, and so $\|a - \mu(a)e\| 
\le kL(a)$.  It
follows that $\|a\|^{\sim} \le kL(a)$, so that the $\rho_L$-radius of 
$S(A)$ is no
greater than $k$.

Suppose conversely that $\|a\|^{\sim} \le rL(a)$ for all $a$.  Let $a 
\in A$ with
$\mu(a) = 0$.  There is a $t \in {\mathbb R}$ such that $\|a - te\| 
\le rL(a)$.  Then
\[
|t| = |\mu(a)-t| = |\mu(a-te)| \le \|a-te\| \le rL(a).
\]
Thus
\[
\|a\| \le \|a-te\| + \|te\| \le 2rL(a).
\]
So for $k = 2r$ we have $\|a\| \le kL(a)$ if $\mu(a) = 0$.
\end{proof}

We see that the constant $k$ is not precisely related to the radius. 
But for our
twisted group algebras there is a very natural state to use, namely the tracial
state $\t$ defined by $\t(f) = f(e)$, which is the vector state for $\d_e \in
\bell^2(G)$.

Suppose now that $\rho_{\o}$ gives $S(C_r^*(G,c))$ finite radius, so 
that as above,
if $\t(f) = 0$ then $\|\pi(f)\| \le 2rL(f)$.  Let $x \in G$ with $x 
\ne e$.  Then
$\t(\d_x) = 0$, and so
\[
1 = \|\pi(\d_x)\| \le 2rL^{\o}(\d_x) = 2r\ell^{\o}(x).
\]
We thus obtain:

\setcounter{proposition}{4}
\begin{proposition}
\label{prop2.5}
If $\rho_{\o}$ gives $S(C_r^*(G,c))$ finite radius $r$, 
then $\ell^{\o}(x) \ge
(2r)^{-1}$ for all $x \ne e$.
\end{proposition}

Thus, for example, if $\th$ is an irrational number, then neither the 
(unbounded)
length function $\ell$ defined on ${\mathbb Z}^2$ by $\ell(m,n) = 
|m+n\th|$, nor
the homomorphism $\o(m,n) = m+n\th$, will give 
metrics for which $S(C^*({\mathbb
Z}^2))$ has finite radius.

But the condition of Proposition \ref{prop2.5} is not at all 
sufficient for finite
radius.  For example, for any $G$ we can define a length function 
$\ell$ by $\ell(x)
= 1$ if $x \ne e$.  Then it is easily checked that if $f = f^*$ then
\[
L^{\ell}(f) = \|f-\t(f)\d_e\|_2.
\]
If $L^{\ell}$ gives $S(C^*(G))$ finite radius, so that there is a 
constant $k$ such
that $\|\pi_f\| \le kL^{\ell}(f)$ if $f(e) = 0$, then it follows that 
$\|\pi_f\| \le
2k\|f\|_2$ when $f(e) = 0$.  Since for any $f$ we have $|f(e)| \le 
\|f\|_2$, it
follows that $\|\pi_f\| \le (2k+1)\|f\|_2$, so that for any $g \in 
C_c(G)$ we have
\[
\|f*g\|_2 \le (2k+1)\|f\|_2\|g\|_2.
\]
This quickly says that the norm on $\bell^2(G)$ can be normalized so that
$\bell^2(G)$ forms an $H^*$-algebra, as defined in section~27 of 
\cite{Lom}.  But
our algebra is unital, and the theory of $H^*$-algebra in \cite{Lom} 
shows that $G$
must then have finite-dimensional square-integrable unitary 
representations.  But
Weil pointed out on page~70 of \cite{Wei} that this means that $G$ is 
compact (so
finite), because if $x \rightarrow U_x$ is the unitary matrix 
representation for a
finite-dimensional square integrable representation, then the matrix 
coefficients of
\[
x \mapsto I = U_xU_x^*
\]
are integrable.

But beyond these elementary comments it is not clear to me what 
happens even for
word-length functions.  Thus we have the basic:

\setcounter{question}{5}
\begin{question}
\label{quest2.6}
For which finitely generated groups $G$ with cocycle $c$ does the word-length
function $\ell$ corresponding to a finite generating subset give a metric
$\rho_{\ell}$ which gives $S(C_r^*(G,c))$ finite diameter?  That is, 
when is there a
constant, $k$, such that if $f \in C_c(G)$ and $f(e) = 0$ then
\[
\|\pi(f)\| \le k\|[M_{\ell},\pi(f)]\|?
\]
(Is the answer independent of the choice of the generating set?)
\end{question}

I do not know the answer to this question when the cocycle $c$ is 
trivial and, for
example, $G$ is the discrete Heisenberg group, or the free group on 
two generators.
In later sections we will obtain some positive answers for $G = 
{\mathbb Z}^d$, but
even that case does not seem easy.

Even less do I know answers to the basic:

\setcounter{question}{6}
\begin{question}
\label{quest2.7}
For which finitely generated groups $G$ with $2$-cocycle $c$ does the 
word-length
function $\ell$ corresponding to a finite generating subset give a metric
$\rho_{\ell}$ which gives $S(C_r^*(G,c))$ the weak-$*$ topology. 
That is \cite{R4},
given that $\rho_{\ell}$ does give $S(C_r^*(G,c))$ finite diameter, when is
\[
{\mathcal B}_1 =
\{f \in C_c(G): \|\pi_f\| \le 1 \mbox{ and } L_{\ell}(f) \le 1\}
\]
a totally-bounded subset of $C_r^*(G)$?
\end{question}

But we now make some elementary observations about this second question.

\setcounter{proposition}{7}
\begin{proposition}
\label{prop2.8}
Let $L$ be a Lip-norm on an order-unit space $A$.  If $L$ is 
continuous for the norm
on $A$, then $A$ is finite-dimensional.
\end{proposition}

\begin{proof}
Much as just above we set
\[
{\mathcal B}_1 =
\{a \in A: \|\a\| \le 1 \mbox{ and } L(a) \le 1\}.
\]
Since $L$ is a Lip-norm, ${\mathcal B}_1$ is totally bounded by
theorem 1.9 of \cite{R4}. But if $L$ is also norm-continuous,
then there is a constant $k \ge 1$ such that $L(a) \le k\|a\|$ for
all $a \in A$. Consequently 
$\{a: \|a\| \le k^{-1}\} \subseteq {\mathcal B}_1$. It follows that
the unit ball for the norm is totally bounded, and so the unit ball 
in the completion of $A$ is compact. But it is well-known that 
the unit ball in a
Banach space is not norm-compact unless 
the Banach space is finite-dimensional.
\end{proof}

\setcounter{corollary}{8}
\begin{corollary}
\label{cor2.9}
Let $A$ be an order-unit space which is represented faithfully as 
operators on a
Hilbert space ${\mathcal H}$.  Let $D$ be a self-adjoint operator on 
${\mathcal H}$,
and set $L(a) = \|[D,a]\|$.  Assume that $L$ is (finite and) a 
Lip-norm on $A$.  If
$D$ is a bounded operator, then $A$ is finite-dimensional.
\end{corollary}

 From this we see that in our setting of $D = M_{\o}$ for 
$C_r^*(G,c)$, if we want
$L^{\o}$ to be a Lip-norm, then we must use unbounded $\o$'s unless 
$G$ is finite.
But it is not clear to me whether $\o$ must always be a proper 
function, that is,
whether $\{x: |\o(x)| \le k\}$ must be finite for every $k$.
However, the referee has pointed out to me that if $\o$ is
actually a length function, then $\o$ must be proper if $L^{\o}$
is to be a Lip-norm. For if it is not proper, then there is a
constant, $r$, with $0 < r \le 1$, such that $S = \{x:\o(x) \le r^{-1}\}$
is infinite. But if $\o$ is a length function then $L(\d_x) = \o(x)$.
Thus $\{r\d_x: x \in S\}$ is a norm-discrete subset of
\[
{\mathcal B}_1 = \{f \in C_c(G): \|f\|_A \le 1 \mbox{ and } L(f) \le 1\},
\]
so that ${\mathcal B}_1$ can not be totally bounded. (See the first 
three sentences of the proof of Lemma 1.4.)

Finally, we will examine briefly three of Connes' axioms for a 
non-commutative Riemannian
geometry \cite{Cn3}.  We begin first with the axiom of ``reality'' 
(axiom $7'$ on
page~163 of \cite{Cn3} and condition~4 on page~483 of \cite{GVF}).  For any
$C^*$-algebra $A$ with trace $\t$ there is a natural and well-known
``charge-conjugation'' operator, $J$, on the GNS Hilbert space for $\t$,
determined by $Ja = a^*$.  We are in that setting, and so our $J$ is given by
\[
(J\xi)(x) = {\bar \xi}(x^{-1})
\]
for $\xi \in \bell^2(G)$.  For any $f \in C_c(G)$ one checks that 
$J\pi_fJ$ is the
operator of {\em right}-convolution by $f^*$, where $f^*(x) = {\bar 
f}(x^{-1})$.  In
particular, $J\pi_fJ$ will commute with any $\pi_g$ for $g \in 
C_c(G)$.  This means
exactly that the axiom of reality is true if 
one considers our geometry to have
dimension $0$.

With the axiom of reality in place, 
Connes requires that $D$ be a ``first-order
operator'' (axiom $2'$ of \cite{Cn3}, 
or condition~5 on page~484 of \cite{GVF},
where the terminology ``first order'' is used).  This axiom requires 
that $[D,a]$
commutes with $JbJ$ for all $a,b \in A$.  For our situation, let 
$\rho_z$ denote
right $c$-twisted translation on $\bell^2(G)$ by $z \in G$, 
so that $J\pi_z^*J =
\rho_z$.  Then in terms of the notation we have established, the first-order
condition requires that $\rho_z$ commutes with $M_{\var_y}$ for each 
$z$ and $y$.
This implies that for each $x \in G$ we have
\[
\o(x) - \o(y^{-1}x) = \o(xz) - \o(y^{-1}xz).
\]
If we choose $z = x^{-1}$ and rearrange, we obtain
\[
\o(x) + \o(y^{-1}) = \o(y^{-1}x) + \o(e).
\]
This says that if we subtract the constant function $\o(e)$, then 
$\o$ is a group
homomorphism from $G$ into ${\mathbb R}$.  Thus the first-order 
condition is rarely
satisfied in our context.  In fact, if we want $\o$ to give 
$S(C_r^*(G))$ finite
radius then it follows from Proposition \ref{prop2.5} that $G \cong 
{\mathbb Z}$ or is finite.

Lastly, we consider the axiom of smoothness (axiom $3$ on page~159 of 
\cite{Cn3}, or
condition~2 on page~482 of \cite{GVF}, where it is called 
``regularity'' rather than
``smoothness'').  This requires that $a$ and $[D,a]$ are in the domains of all
powers of the derivation $T \mapsto [|D|,T]$.  In our context $|D| = 
M_{|\o|}$.  But
\[
||\o(x)| - |\o(z^{-1}x)|| \le |\o(x)-\o(z^{-1}x)|,
\]
so that $|\o|$ is translation-bounded when $\o$ is.  From this it is 
easily seen
that the axiom of smoothness is always satisfied in our setting.

%%%%%%%%%%%%%%%%%%%%%%%%%%%%%%%%%%%%%%%%%
\setcounter{section}{2}
\section{The cosphere algebra}
\label{sec3}
%%%%%%%%%%%%%%%%%%%%%%%%%%%%%%%%%%%%%%%%%%%

We now begin to establish some constructions which will permit us to 
obtain positive
answers to Questions \ref{quest2.6} and \ref{quest2.7} 
for the groups ${\mathbb
Z}^d$, and which may eventually be helpful in dealing with other groups.

Connes has shown (section~6 of \cite{C4}, \cite{GLc}) how to 
construct for each
spectral triple $(A,{\mathcal H},D)$ a certain $C^*$-algebra, denoted 
$S^*A$.  He
shows that if $A = C^\infty({\mathcal M})$ 
where ${\mathcal M}$ is a compact 
Riemannian spin
manifold, and if $({\mathcal H},D)$ is the corresponding Dirac 
operator, then $S^*A$
is canonically isomorphic to the algebra of continuous functions on the unit
cosphere bundle of ${\mathcal M}$.  Thus in the general case it seems 
reasonable to
call $S^*A$ the cosphere algebra of $(A,{\mathcal H},D)$.  (In 
\cite{GLc} $S^*A$ is
called the ``unitary cotangent bundle''.)  In this section we will 
explore what this
cosphere algebra is for our (almost) spectral triples of form
$(C_c(G),\bell^2(G),M_{\o})$.  (I thank Pierre Julg for helpful 
comments about this
at an early stage of this project.)

We now review the general construction.  But for our purposes we do 
not need the
usual further hypothesis of finite summability for $D$.  Thus we just 
require that
we have $(A,{\mathcal H},D)$ such that $[D,a]$ is bounded for all $a 
\in A$.  But,
following Connes, we also make the smoothness requirement that $[|D|,a]$ be 
bounded for all $a
\in A$.  We saw in the previous section that this latter condition is always
satisfied in our setting where $D = M_{\o}$.

Connes' construction of the algebra $S^*A$ is as follows.  (See also the
introduction of \cite{GLc}.)  Form the strongly continuous 
one-parameter unitary
group $U_t = \exp(it|D|)$.  Let ${\mathcal C}_D$ be the $C^*$-algebra 
of operators
on ${\mathcal H}$ generated by the algebra ${\mathcal K}$ of compact 
operators on
${\mathcal H}$ together with all of the algebras $U_tAU_{-t}$ for $t 
\in {\mathbb
R}$.  (Note that usually $U_tAU_{-t} \not\subseteq A$.)  Clearly the action of
conjugating by $U_t$ carries ${\mathcal C}_D$ into itself.  We denote 
this action of
${\mathbb R}$ on ${\mathcal C}_D$ by $\eta$.  Because of the requirement that
$[|D|,a]$ be bounded, the action $\eta$ is strongly continuous on 
${\mathcal C}_D$. (See the first line of the proof of corollary 10.16
of \cite{GVF}.)
Since ${\mathcal K}$ is an ideal ($\eta$-invariant) in ${\mathcal 
C}_D$, we can form
${\mathcal C}_D/{\mathcal K}$.  Then by definition $S^*A = {\mathcal 
C}_D/{\mathcal
K}$.  The action $\eta$ drops to an action of ${\mathbb R}$ on 
$S^*A$, which Connes
calls the ``geodesic flow''.

We now work out what the above says for our case in which we have
$(C_r^*(G,c),\bell^2(G),M_{\o})$.  
We will write ${\mathcal C}_{\o}$ instead of
${\mathcal C}_D$.  Since only $|\o|$ is pertinent, we assume for a 
while that $\o
\ge 0$.  Set $u_t(x) = \exp(it\o(x))$ for $t \in {\mathbb R}$, so 
that the $U_t$ of
the above construction becomes $M_{u_t}$.  Then for each $y \in G$ our algebra
${\mathcal C}_{\o}$, defined as above, must contain
\[
U_t\pi_yU_t^* = M_{u_t}M_{\a_y(u_t^*)}\pi_y = M_{u_t\a_y(u_t^*)} \pi_y.
\]
But ${\mathcal C}_{\o}$ must also contain $(\pi_y)^{-1}$, and thus it 
contains each
$u_t\a_y(u_t^*)$, where for notational simplicity we omit $M$.  But
\[
(u_t\a_y(u_t^*))(x) = \exp(it(\o(x) - \o(y^{-1}x))) = \exp(it\var_y(x)).
\]
Since $\var_y$ is bounded, the derivative of $U_t\a_y(U_t^*)$ at $t = 
0$ will be
the norm-limit of the difference quotients.  Thus we see that also $\var_y \in
{\mathcal C}_{\o}$ for each $y \in G$.  But ${\mathcal C}_{\o} 
\supseteq {\mathcal
K}$, and so ${\mathcal C}_{\o} \supseteq C_{\infty}(G)$, the space
of continuous functions vanishing at infinity, where the elements of
$C_{\infty}(G)$ are here viewed as multiplication operators.  Note also that
${\mathcal C}_{\o}$ contains the identity element.

All of this suggests that we consider, inside the algebra $C_b(G)$ of bounded
functions on $G$, the unital norm-closed subalgebra generated by 
$C_{\infty}(G)$
together with $\{\var_y: y \in G\}$.  We denote this subalgebra by 
$E_{\o}$.  Let
${\bar G}^{\o}$ denote the maximal ideal space of $E_{\o}$, with its compact
topology, so that $E_{\o} = C({\bar G}^{\o})$.  Note that $G$ sits in ${\bar
G}^{\o}$ as a dense open subset because $E_{\o} \supseteq 
C_{\infty}(G)$.  That is,
${\bar G}^{\o}$ is a compactification of the discrete set $G$.  We 
will call it the
$\o$-{\em compactification} of $G$.  Note that $C({\bar G}^{\o})$ is separable
because $G$ is countable and so there is 
only a countable number of $\var_y$'s.
Thus the compact topology of ${\bar G}^{\o}$ has a countable base.

The action $\a$ of $G$ on $C_b(G)$ by left translation 
clearly carries $E_{\o}$
into itself.  From this we obtain an induced action on ${\bar G}^{\o}$ by
homeomorphisms.  We denote this action again by $\a$.

Of course $C({\bar G}^{\o})$ is faithfully represented as an algebra 
of pointwise
multiplication operators on $\bell^2(G)$.  This representation, $M$, 
together with
the representation $\pi$ of $G$ on $\bell^2(G)$ form a covariant representation
\cite{Pdr}, \cite{ZM} of $(C({\bar G}^{\o}),G,\a,c)$.  We have already seen
earlier several instances of the covariance relation $\pi_xM_f =
M_{\a_x(f)}\pi_x$.  
The integrated form of this covariant representation, which
we denote again by $\pi$, gives then a representation on $\bell^2(G)$ 
of the full
twisted crossed product algebra $C^*(G,C({\bar G}^{\o}),\a,c)$.  It is clear
from the above discussion that our algebra ${\mathcal C}_{\o}$ contains
$\pi(C^*(G,C({\bar G}^{\o}),\a,c))$.  
But for any $y \in G$ and $t \in {\mathbb
R}$ we have $\exp(it\var_y) \in C({\bar G}^{\o})$.  From our earlier 
calculation
this means that $\pi(C^*(G,C({\bar G}^{\o}),\a,c))$ contains $U_t\pi_yU_t^*$.
Thus it also contains $U_t\pi(C_c(G))U_t^*$.  Consequently:

\setcounter{lemma}{0}
\begin{lemma}
\label{lem3.1}
We have ${\mathcal C}_{\o} = \pi(C^*(G,C({\bar
G}^{\o}),\a,c))$.
\end{lemma}

Now $C({\bar G}^{\o})$ contains $C_{\infty}(G)$ as an $\a$-invariant 
ideal.  The
following fact must be known, but I have not found a reference for it.

\setcounter{lemma}{1}
\begin{lemma}
\label{lem3.2}
With notation as above,
\[
C^*(G,C_{\infty}(G),\a,c) \cong {\mathcal K}(\bell^2(G)),
\]
the algebra of compact operator on $\bell^2(G)$, 
with the isomorphism given by $\pi$.
\end{lemma}

\begin{proof}
If we view elements of $C_c(G,C_{\infty}(G))$ as functions on $G \x G$, and if
for $f \in C_c(G,C_{\infty}(G))$ we set $(\Phi f)(x, y) = f(x,y)c(x,x^{-1}y)$,
then
\[
\Phi(f*_cg) = (\Phi f)*(\Phi g),
\]
where only here we let $*_c$ denote convolution (in the crossed 
product) twisted by
$c$, while $*$ denotes ordinary convolution.  The verification 
requires using the
$2$-cocycle identity to see that
\[
c(y,y^{-1}z)c(y^{-1}x,x^{-1}z) = c(y,y^{-1}x)c(x,x^{-1}z).
\]
The untwisted crossed product $C_{\infty}(G)\times_{\a}G$ is well-known to 
be carried
onto ${\mathcal K}(\bell^2(G))$ by $\pi$.  (See \cite{R12}.)  
(For non-discrete
groups one must be more careful, because cocycles are often only 
measurable, not
continuous.)
\end{proof}

Because ${\mathcal K}(\bell^2(G))$ is simple, it follows that the reduced
$C^*$-algebra $C_r^*(G,C_{\infty}(G),\a,c)$ coincides with the full twisted
crossed product, even when $G$ is not amenable.  Anyway, the 
consequence of this
discussion is:

\setcounter{proposition}{2}
\begin{proposition}
\label{prop3.3}
With notation as above, the cosphere algebra is
\[
S^*A = \pi(C^*(G,C({\bar G}^{\o}),\a,c))/{\mathcal K}(\bell^2(G)).
\]
\end{proposition}

For an element of $\pi(C^*(G,C({\bar G}^{\o}),\a,c)$ 
it is probably appropriate
to call its image in $S^*A$ its ``symbol'', in analogy with the situation for
pseudodifferential operators.

We can use recently-developed technology to obtain a simpler picture 
in those cases
in which the action $\a$ of $G$ on ${\bar G}^{\o}$ is amenable \cite{AD1},
\cite{ARn}, \cite{AD2}, \cite{HgR}, \cite{Hgs}.  This action will 
always be amenable
if $G$ itself is amenable, which will be the case when we consider 
${\mathbb Z}^d$
in detail later.  So the following comments will only be needed there for that
case.  But we will see in Section~10 that the action can be amenable 
also in some
situations for which $G$ is not amenable, namely for the free group on two
generators and its standard word-length function.

Let $\p_{\o}G = {\bar G}^{\o} \setminus G$.  It is reasonable to call
$\p_{\o}G$ the ``$\o$-boundary'' of $G$.  Notice that $\a$ carries
$\p_{\o}G$ into itself.  Suppose that the action $\a$ of $G$ on
$\p_{\o}G$ is amenable \cite{AD2}, \cite{ARn}.  One of the equivalent 
conditions for
amenability of $\a$ (for discrete $G$) is that the quotient map from
$C^*(G,C(\p_{\o}G))$ onto $C_r^*(G,C(\p_{\o}G))$ is an isomorphism 
(theorem $4.8$ of
\cite{AD1} or theorem $3.4$ of \cite{AD2}).  (No cocycle $c$ is 
involved here.)  In
proposition $2.4$ of \cite{KmY} it is shown that for situations like this
the amenability of the action on $\p_{\o}G$ is equivalent to amenability 
of the action
on ${\bar G}^{\o}$.  (I thank Claire Anantharaman--Delaroche for bringing this
reference to my attention, and I thank both her and Jean Renault for helpful
comments on related matters.)  The proof in \cite{KmY} uses the 
characterization of
amenability of the action in terms of nuclearity of the crossed 
product.  Here is
another argument which does not use nuclearity.  Following remark $4.10$ of
\cite{Rn}, we consider the exact sequence of full crossed products
\[
0 \rightarrow C^*(G,C_{\infty}(G),\a) \rightarrow C^*(G,C({\bar G}^{\o}),\a)
\rightarrow C^*(G, C(\p_{\o}G),\a) \rightarrow 0
\]
and its surjective maps onto the corresponding sequence of reduced 
crossed products
(which initially is not known to be exact).  A simple diagram-chase 
shows that if
the quotient map onto $C_r^*(G,C(\p_{\o}G),\a)$ is in fact an 
isomorphism, then the
sequence of reduced crossed products is in fact exact.  Also, as 
discussed above,
$C^*(G,C_{\infty}(G),\a)$ is the algebra of compact operators, so 
simple, and so the
quotient map from it must be an isomorphism.  A second simple 
diagram-chase then
shows that the quotient map from $C^*(G,C({\bar G}^{\o}),\a)$ must be an
isomorphism, so that the action $\a$ of $G$ on ${\bar G}^{\o}$ is 
amenable.  (The
verification that if the action on ${\bar G}^{\o}$ is amenable then 
so is that  on
$\p_{\o}G$ follows swiftly from the equivalent definition of 
amenability in terms of
maps whose values are probability measures on $G$.  This definition 
is given further
below and in example $2.2.14(2)$ of \cite{ARn}.)

For our general functions $\o$ it is probably not reasonable to hope 
to find a nice
criterion for amenability of the action.  But in the case in which $\o$ is a
length-function $\ell$ (in which case we write $\p_{\ell}G$ instead 
of $\p_{\o}G$),
we will obtain in the next sections considerable information about 
$\p_{\ell}G$, and
so it is reasonable to pose:

\setcounter{question}{3}
\begin{question}
\label{quest3.4}
Let $G$ be a finitely generated group, and let $\ell$ be the 
word-length function
for some finite set of generators.  Under what conditions will the 
action of $G$ on
$\p_{\ell}G$ be amenable?  For which class of groups will there exist 
a finite set
of generators for which the action is amenable?  For which class of 
groups will this
amenability be independent of the choice of generators?
\end{question}

It is known that if $G$ is a word-hyperbolic group, then its action 
on its Gromov
boundary is amenable.  See the appendix of \cite{ARn}, written by 
E.~Germain, and
the references given there.  We would have a positive answer to Question
\ref{quest3.4} for word-hyperbolic groups if we had a positive answer to:

\setcounter{question}{4}
\begin{question}
\label{quest3.5}
Is it the case that for any word-hyperbolic group $G$ and any 
word-length function
on $G$ for a finite generating set, there is an equivariant 
continuous surjection
from $\p_{\ell}G$ onto the Gromov boundary of $G$?
\end{question}

This seems plausible in view of our discussion of geodesic rays in 
the next section,
since the Gromov boundary considers geodesic rays which stay a finite 
distance from
each other to be equivalent.

We now explore briefly the consequences of the action being amenable. 
The first
consequence is that the full and reduced twisted crossed products 
coincide.  We have
discussed the case of a trivial cocycle $c$ above.
I have not seen the twisted case stated in the literature, but it 
follows easily
from what is now known.  We outline the proof.  To every $2$-cocycle there is
associated an extension, $E$, of $G$ by ${\mathbb T}$.  As a 
topological space $E =
{\mathbb T} \x G$, and the product is given by $(s,x)(t,y) = 
(stc(x,y),xy)$.  (See
III.5.12 of \cite{FD}.)  We can compose the evident map from $E$ onto 
$G$ with $\a$
to obtain an action, $\a$, of $E$ on ${\bar G}^{\o}$.  Let $W$ be any 
compact space
on which $G$ acts, with the action denoted by $\a$.  If $\a$ is 
amenable, then by
definition (example $2.2.14(2)$ of \cite{ARn}, \cite{AD2}, 
\cite{HgR}, \cite{Hgs})
there is a sequence $\{m_j\}$ of weak-$*$ continuous maps from $W$ 
into the space of
probability measures on $G$ such that, 
for $\a$ denoting also the corresponding
action on probability measures, we have for every $x \in G$
\[
\lim_j \sup_{w \in W} \|\a_x(m_j(w)) - m_j(\a_x(w))\|_1 = 0.
\]
Let $h$ denote normalized Haar measure on ${\mathbb T}$, and for each 
$j$ and each
$w \in W$ let $n_j(w)$ be the product measure $h \otimes m_j(w)$ on 
$E$.  Thus each
$n_j(w)$ is a probability measure on $E$.  It is easily verified that 
the function
$w \mapsto n_j(w)$ is weak-$*$ continuous.  Furthermore, a straight-forward
calculation shows that
\[
\a_{(s,x)}(n_j) = h \otimes \a_x(m_j)
\]
for each $(s,x) \in E$ and each $j$.  Now $E$ is not discrete.  But from this
calculation it is easily seen that the action of $E$ on $W$ is 
amenable, where now
we use definition $2.1$ of \cite{AD2}.  Then from theorem $3.4$ of \cite{AD2}
(which is a special case of proposition $6.1.8$ of \cite{ARn}), 
it follows that
$C^*(E,C(W),\a)$ coincides with $C_r^*(E,C(W),\a)$.

Now let $p$ be the function on ${\mathbb T}$ 
defined by $p(t) = \exp(2\pi it)$,
where here we identify ${\mathbb T}$ with ${\mathbb R}/{\mathbb Z}$.  Since
${\mathbb T}$ is an open subgroup of $E$, we can view $p$ as a 
function on $E$ by
giving it value $0$ off of ${\mathbb T}$.  Since ${\mathbb T}$ is 
central in $E$,
and $\a$ is trivial on ${\mathbb T}$, and $C(W)$ is unital, it 
follows that $p$ is a
central projection in $C^*(E,C(W),\a)$.  
From this it follows that the cut-down
algebras $pC^*(E,C(W),\a)$ and $pC_r^*(E,C(W),\a)$ coincide.  But it 
is easily seen
(see page~84 of \cite{Fl1} or page~144 of \cite{Fl2}) that $pC^*(E,C(W),\a) =
C^*(G,C(W),\a,c)$, and similarly for $C_r^*$.  In this way we obtain:

\setcounter{proposition}{5}
\begin{proposition}
\label{prop3.6}
Let $G$ be a discrete group, let $\a$ be an action of $G$ on a 
compact space $W$,
and let $c$ be a $2$-cocycle on $G$.  If the action $\a$ is amenable, then
$C^*(G,C(W),\a,c)$ coincides with $C_r^*(G,C(W),\a,c)$.
\end{proposition}

With some additional care the above proposition can be extended to 
the case in which
$W$ is only locally compact.  In that case the projection $p$ is only in the
multiplier algebras of the twisted crossed products.

We now return to the case in which $G$ acts on ${\bar G}^{\o}$ and 
$\p_{\o}G$.  From
the above proposition it follows that if $G$ acts amenably on 
$\p_{\o}G$, and so on
${\bar G}^{\o}$, then we can view $\pi$ as a representation of the 
reduced crossed
product $C_r^*(G,C({\bar G}^{\o}),\a,c)$.  This has the benefit that 
we can apply
corollary $4.19$ of \cite{ZM} to conclude that $\pi$ is a faithful 
representation of
$C_r^*(G,C({\bar G}^{\o}),\a,c)$.  The hypotheses of this corollary 
$4.19$ are that
$M$ be a faithful representation of $C({\bar G}^{\o})$, which is 
clearly true, and
that $M$ be $G$-almost free (definition $1.12$ of \cite{ZM}).  This 
latter means
that for any non-zero subrepresentation $N$ of $M$ and any $x \in G$ 
with $x \ne e$
there is a non-zero subrepresentation $P$ of $N$ whose composition 
with the inner
automorphism from $x$ is disjoint from $P$.  But subrepresentations of $M$
correspond to non-empty subsets of $G$, 
and for $P$ we can take any one-point subset
of a given subset.  Thus our algebra ${\mathcal C}^{\o}$ coincides 
(under $\pi$)
with $C^*(G,C({\bar G}^{\o}),\a,c)$.

Now from Lemma \ref{lem3.2} we know that $C^*(G,C_{\infty}(G),\a,c)$ 
coincides with
${\mathcal K}(\bell^2(G))$, and the process of forming full twisted 
crossed products
preserves short exact sequences.  (See the top of page~149 of 
\cite{ZM}.)  Thus
from Proposition \ref{prop3.3}, and on removing our requirement that 
$\o \ge 0$, we
obtain:

\setcounter{theorem}{6}
\begin{theorem}
\label{th3.7}
Let $\o$ be a translation bounded function on $G$ such that the 
action of $G$ on
$\p_{|\o|}G$ is amenable.  Then the cosphere algebra $S_{\o}^*A$ for
$(C_r^*(G,c),\bell^2(G),M_{\o})$ is (naturally identified with)
\[
S_{\o}^*A = C^*(G,C(\p_{\o}G),\a,c) = C_r^*(G,C(\p_{\o}G),\a,c).
\]
\end{theorem}

%%%%%%%%%%%%%%%%%%%%%%%%%%%%%%%%%%%%%%%%%%%%%%
\setcounter{section}{3}
\section{The metric compactification}
\label{sec4}
%%%%%%%%%%%%%%%%%%%%%%%%%%%%%%%%%%%%%%%%%%%%%%%%

The purpose of this section is to show that when $\o$ is a 
length-function on $G$
then geodesic rays in $G$ for the metric on $G$ from $\o$ give points in the
compactification ${\bar G}^{\o}$.  This will be a crucial tool for us 
in dealing
with ${\mathbb Z}^d$, since it will supply us with a sufficient collection
of points in the boundary which have finite orbits.  
We will also see that ${\bar G}^{\o}$ is then a 
special case of a
compactification of complete locally compact 
metric spaces introduced by Gromov
\cite{G3} some time ago.  
(This is probably related to the comment which Connes
makes about nilpotent groups in the second paragraph after the end of 
the proof of
proposition~2 of section~6 of \cite{C5}.)  Gromov's definition appears fairly
different from that which we gave in the previous section, and so our 
treatment here
can also be viewed as showing how to define Gromov's compactification 
as the maximal
ideal space of a unital commutative $C^*$-algebra.  We will refrain 
from using here
the terms ``Gromov compactification'' and ``Gromov boundary'', since 
these terms
seem already reserved in the literature for use with hyperbolic 
spaces, where they
have a different meaning and give objects which depend only on the coarse
quasi-isometry class of the metric.  (See IIIH3 of \cite{BrH}.)  We 
will instead use
the terms ``metric compactification'' and ``metric boundary'', and 
our notation will
often show the dependence on the metric.  We will see in Example 
\ref{exam5.2} that
for a hyperbolic metric space the metric boundary and the Gromov 
boundary can fail
to be homeomorphic.

Let $(X,\rho)$ be a metric space, and let $C_b(X)$ denote the algebra 
of continuous
bounded functions on $X$, equipped with the supremum norm $\|\cdot\|_{\infty}$.
Motivated by the observations in the previous section, we define 
$\var_{y,z}$ on $X$
for $y,z \in X$ by
\[
\var_{y,z}(x) = \rho(x,y) - \rho(x,z).
\]
Then the triangle inequality tells us that $\|\var_{y,z}\|_{\infty} 
\le \rho(y,z)$,
so that $\var_{y,z} \in C_b(X)$.  But on setting $x = z$ we see that, in fact,
$\|\var_{y,z}\|_{\infty} = \rho(y,z)$.  Let $H_{\rho}$ denote the 
linear span in
$C_b(X)$ of $\{\var_{y,z}: y,z \in X\}$.  Suppose that we fix some 
base point $z_0
\in X$.  Then it is easily checked that $\var_{y,z} = \var_{z_0,z} - 
\var_{z_0,y}$.
Thus $H_{\rho}$ is equally well the linear span of $\{\var_{z_0,y}: y 
\in X\}$, but
is independent of the choice of $z_0$.  (It will be useful to us that 
we can change
base-points at will.)  We often find it convenient to fix $z_0$, and to set
$\var_y = \var_{z_0,y}$, so that $H_{\rho}$ is the linear span of the
$\var_y$'s.  When $X$ is a group, it is natural to choose $z_0 = e$.  We were
implicitly doing this in the previous section.  We note that 
$\|\var_y\|_{\infty} =
\rho(y,z_0)$.

Much as above, we have $\var_y - \var_z = \var_{z,y}$, and so $\|\var_y -
\var_z\|_{\infty} = \|\var_{z,y}\|_{\infty} = \rho(y,z)$.  Thus the mapping $y
\mapsto \var_y$ is an isometry from $(X,\rho)$ into $C_b(X)$.  The 
latter space is
complete, and so this isometry extends to the completion of $X$.

We desire to obtain a compactification of $X$ to which all of the
functions $\var_y$ extend as continuous functions. We want $X$
to be an open subset of the compactification, and so 
we must require that $X$ is locally compact. Then the 
various compactifications
of $X$ in which $X$ is open are just the maximal-ideal spaces
of the various unital closed *-subalgebras of $C_b(X)$
which contain $C_\infty(X)$. Thus we set:

\setcounter{definition}{0}
\begin{definition}
\label{def4.1}
Let $(X,\rho)$ be a metric space whose topology is locally compact. 
Let ${\mathcal
G}(X,\rho)$ be the norm-closed subalgebra of $C_b(X)$ which is generated by
$C_{\infty}(X)$, the constant functions, and $H_{\rho}$.  
Let ${\bar X}^{\rho}$
denote the maximal ideal space of ${\mathcal G}(X,\rho)$.  We call 
${\bar X}^{\rho}$
the {\em metric compactification} of $X$ for $\rho$.
\end{definition}

Then, essentially by construction, ${\bar X}^{\rho}$ is a 
compactification of $X$ (within which $X$ is open). We remark
that if, instead, we take the norm-closed subalgebra of $C_b(X)$
generated by all of the bounded Lipschitz functions, then we
obtain the algebra of all bounded uniformly continuous (for $\rho$)
functions on $X$. (See the bottom of page 23 of \cite{Wea}.)

 It is natural to think of 
$X_{\rho}{\backslash}X$ as a
boundary at infinity for $X$.  But from a metric standpoint this is not always
reasonable.  Suppose that $X$ is not complete.  Each of 
the functions
$\var_y$ is a Lipschitz function, and so extends to the completion 
${\hat X}^{\rho}$
of $X$.  Each $f \in C_{\infty}(X)$ extends continuously to ${\hat 
X}^{\rho}$ by
setting it equal to $0$ off $X$.  The constant functions obviously 
extend to ${\hat
X}^{\rho}$.  Thus the algebraic algebra generated by $H_{\rho}$, 
$C_{\infty}(X)$ and
the constant functions extends to an algebra of functions on ${\hat 
X}^{\rho}$, and
the supremum norm is preserved under this extension.  Thus our 
completed algebra
${\mathcal G}(X,\rho)$ can be viewed as a unital subalgebra of $C_b({\hat
X}^{\rho})$.  It is easily seen that this algebra separates the 
points of ${\hat
X}^{\rho}$.  (E.g., use the fact that $\rho$ extends to the 
completion.)  Thus we
obtain a (continuous) injection of ${\hat X}^{\rho}$ into ${\bar 
X}^{\rho}$.  But
there is no reason that ${\hat X}^{\rho}$ should be open in ${\bar X}^{\rho}$,
notably if the completion is not locally compact.  
Even if ${\hat X}^{\rho}$ is
locally compact, the points of ${\hat X}^{\rho}{\backslash}X$ will 
all be of finite
distance from the points of $X$, and so are not ``at infinity''.  For 
this reason it
seems best to define the ``boundary'' only for {\em complete} locally compact
metric spaces.  Thus we make:

\setcounter{definition}{1}
\begin{definition}
\label{def4.2}
Let $(X,\rho)$ be a metric space which is complete and locally 
compact.  Then its
{\em metric boundary} is ${\bar X}^{\rho}{\backslash}X$.  We will 
denote the metric
boundary by $\p_{\rho}X$.
\end{definition}

We now show that the metric compactification and the 
metric boundary which we have
defined above coincide with those 
constructed by Gromov \cite{G3} in a somewhat
different way.  Gromov proceeds as follows.  (See also $3.1$ of 
\cite{BGS}, II.1 of
\cite{Bll} and II.8.12 of \cite{BrH}.)  Let $(X,\rho)$ be a complete 
locally compact
metric space, let $C(X)$ denote the vector space of all continuous (possibly
unbounded) functions on $X$, and equip $X$ with the topology of 
uniform convergence
on compact subsets of $X$.  Let $C_*(X)$ denote the quotient of $C(X)$ by the
subspace of constant functions, with the quotient topology.  For $f 
\in C(X)$ denote
its image in $C_*(X)$ by ${\bar f}$.  For $y \in X$ set $\psi_y(x) = 
\rho(x,y)$.  Then
$x \mapsto \psi_x$ is an embedding of $X$ into $C(X)$.  Let $\iota$ denote the
corresponding embedding of $X$ into $C_*(X)$, and let ${\mathcal 
C}\ell(X)$ be the
closure of $\iota(X)$ in $C_*(X)$.  Then ${\mathcal C}\ell(X)$ can be 
shown to be
compact, and $\iota(X)$ can be shown to be open in ${\mathcal 
C}\ell(X)$, so that
${\mathcal C}\ell(X){\backslash}X$ is a boundary at infinity for $X$.

We now explain the relationship between this construction of Gromov and our
construction given earlier in this section.  Fix a base point $z_0$. 
For any given
$u \in {\bar X}^{\rho}$ define the function $g_u$ by $g_u(x) = 
-\var_x(u)$, where
$\var_x$ is now viewed as a function on ${\bar X}^{\rho}$.  If $u \in X$ then
$g_u(x) = \rho(u,x) - \rho(u,z_0)$.  Since $\rho(u,z_0)$ is constant 
in $x$, the
image of $g_u$ in $C_*(X)$ is exactly Gromov's $\iota(u)$.  On the other hand,
suppose that $u \in \p_{\rho}X$.  Because $X$ is dense in ${\bar 
X}^{\rho}$, there
is a net $\{y_{\a}\}$ of elements of $X$ which converges to $u$.  
Then for each $x \in X$ we have
\[
g_u(x) = -\var_x(u) = -\lim \var_x(y_{\a}) = \lim g_{y_{\a}}(x).
\]
That is, $g_{y_{\a}}$ converges to $g_u$ pointwise on $X$.  But each 
$g_y$ for $y
\in X$ is clearly a Lipschitz function of Lipschitz constant $1$, 
and pointwise
convergence of a net of functions of bounded Lipschitz constant 
implies uniform
convergence on compact sets.  
Thus $g_{y_{\a}}$ converges uniformly to $g_u$ on 
compact subsets
of $X$, so that ${\bar g}_u \in {\mathcal C}\ell(X)$.  (In the 
literature cited above,
$g_u$ would be called a {\em horofunction} if $u \in \p_{\rho}X$.) 
In this way we
obtain a mapping, $u \mapsto {\bar g}_u$, from ${\bar X}^{\rho}$ to ${\mathcal
C}\ell(X)$.  If ${\bar g}_u = {\bar g}_v$ for some $u,v \in {\bar 
X}^{\rho}$, then
there is a constant, $k$, such that $\var_x(u) = \var_x(v) + k$ for 
all $x \in X$.
From this it is easily seen that $u = v$.  Thus the mapping $u 
\mapsto {\bar g}_u$
is injective on ${\bar X}^{\rho}$.  Finally, 
if $\{u_{\a}\}$ is a net in ${\bar
X}^{\rho}$ which converges to $u \in {\bar X}^{\rho}$, then, much as above,
$g_{u_{\a}}$ converges to $g_u$ pointwise, and so uniformly on 
compact sets.  Thus
the mapping $u \mapsto {\bar g}_u$ is continuous from ${\bar X}^{\rho}$ into
${\mathcal C}\ell(X)$.  Since ${\bar X}^{\rho}$ is compact, 
it follows that this mapping is a 
homeomorphism onto its
image.  But the image of $X$ in ${\mathcal C}\ell(X)$ is dense, and 
so the mapping
is a homeomorphism from ${\bar X}^{\rho}$ onto ${\mathcal C}\ell(X)$, 
and so from
$\p_{\rho}X$ to ${\mathcal C}\ell(X){\backslash}X$, as desired.

For our later purposes it is important for us to examine the 
relationship between
geodesics and points of $\p_{\rho}X$.  Much of the content of the 
next paragraphs
appears in some form in various places in the literature \cite{BGS}, 
\cite{Bll},
\cite{BrH}, though usually not in the generality we consider here.  
And here we
reformulate it in terms of our approach to the construction of
$\p_{\rho}X$.

We will not assume that our metric spaces are connected.  For 
example, we will later
consider ${\mathbb Z}^d$ with its Euclidean metric from ${\mathbb 
R}^d$.  Every ray
(half-line) in ${\mathbb R}^d$ should give a direction toward 
infinity for ${\mathbb
Z}^d$.  But if the direction involves irrational angles, the ray may not meet
${\mathbb Z}^d$ at an infinite number of points.  So we need a slight 
generalization of geodesic rays. For perspective we also include a
yet weaker definition.

\setcounter{definition}{2}
\begin{definition}
\label{def4.3}
Let  $(X,\rho)$ be a metric space, let $T$ be an unbounded
subset of ${\mathbb R}^+$ which contains $0$, and let $\g$ be
a function from $T$ into $X$. We will say that:

\begin{itemize}
\item[{\rm a)}] 
$\g$ is a {\em geodesic ray} if 
$\rho(\g(t),\g(s)) = |t-s|$ for all $t, s \in T$.
\item[{\rm b)}] 
$\g$ is an {\em almost-geodesic ray} if it
satisfies the condition:

For every $\e > 0$ there is an integer $N$ such that if $t,s \in T$ and $t
\ge s \ge N$, then
\[
|\rho(\g(t),\g(s)) + \rho(\g(s),\g(0)) - t| < \e.
\]

\item[{\rm c)}]
$\g$ is a {\em weakly-geodesic ray} if for every $y \in X$ and
every $\e > 0$ there is an integer $N$ such that if $s, t \ge N$ then
\[
|\rho(\g(t),\g(0)) - t| < \e
\]
and
\[
|\rho(\g(t), y) - \rho(\g(s), y) - (t - s)| < \e.
\]
\end{itemize}
\end{definition}

It is evident that any geodesic ray is an almost-geodesic ray. (I thank
Simon Wadsley for pointing out to me that my definition of 
weakly-geodesic rays in the first version of this paper was defective.)

\setcounter{lemma}{3}
\begin{lemma}
\label{lem4.4}
Let $\g$ be an almost-geodesic ray, and let $(\e,N)$ be as in Definition
4.3b.  Then for $t \ge s \ge N$ we have:

\begin{itemize}
\item[{\rm a)}] $|\rho(\g(t),\g(0)) - t| < \e$.
\item[{\rm b)}] $|\rho(\g(t),\g(s)) - (t-s)| < 2\e$.
\item[{\rm c)}] $\rho(\g(t),\g(s)) < \rho(\g(t),\g(0)) - 
\rho(\g(s),\g(0)) + 2\e$.
\end{itemize}
\end{lemma}

\begin{proof}
For a) set $s = t$ in the condition of Definition 4.3b.  For b) we have
\begin{eqnarray*}
&& |\rho( \g(t)),\g(s)) - (t-s)| \\
&& \hskip0.2in = |(\rho(\g(t),\g(s)) + \rho(\g(s),\g(0)) - t) 
- (\rho(\g(s),\g(0)) - s)| < 2\e.
\end{eqnarray*}
Finally, for c) we have
\begin{eqnarray*}
&&\rho(\g(t), \g(s)) \\
&& \hskip0.3in = (\rho(\g(t),\g(s)) + \rho(\g(s),\g(0)) - t) 
- \rho(\g(s),\g(0)) \\
&& \hskip0.3in + \rho(\g(t),\g(0)) 
- (\rho(\g(t),\g(0)) - t) \\
&& \hskip0.3in < \rho(\g(t),\g(0)) - \rho(\g(s),\g(0)) + 2\e.
\end{eqnarray*}

\end{proof}

\setcounter{lemma}{4}
\begin{lemma}
\label{lem4.5}
Any almost-geodesic ray is weakly geodesic. 
Let $\g$ be a weakly-geodesic ray. Take $\g(0)$ as the base-point
for defining $\var_y$ for any $y \in X$. Then
$\lim_{t \rightarrow \infty}\var_y(\g(t))$ exists for every
$y \in X$. If $\g$ is
actually a geodesic ray, then $t \mapsto \var_y(\g(t))$ is
a non-decreasing (bounded) function.
\end{lemma}

\begin{proof}
To motivate the rest of the proof, suppose first that $\g$ is
a geodesic ray. We show that $t \mapsto \var_y(\g(t))$ is a
non-decreasing function (so has a limit). For $t \ge s$ we have
\begin{eqnarray*}
\var_y(\g(t)) &- &\var_y(\g(s)) \\
&= &\rho(\g(t), \g(0)) - \rho(\g(t), y) -  
\rho(\g(s), \g(0)) + \rho(\g(s), y)\\
&= & t - s + \rho(\g(s), y) - \rho(\g(t), y)\\
&= &\rho(\g(t), \g(s)) + \rho(\g(s), y) - \rho(\g(t), y) \ge 0
\end{eqnarray*}
by the triangle inequality.

Next, let $\g$ be an almost-geodesic ray. It is useful and
instructive to first see why $\lim_{t \rightarrow \infty}\var_y(\g(t))$
exists.
Given $\e > 0$, take $N$ as in Definition 4.3b. We will show first that
if $t \ge s \ge N$ then $\var_y(\g(t)) > \var_y(\g(s)) - 3\e$.
In fact,
\begin{eqnarray*}
\var_y(\g(t)) &- &\var_y(\g(s)) \\
&= &\rho(\g(t),\g(0)) - \rho(\g(t),y) 
- \rho(\g(s),\g(0)) + \rho(\g(s),y) \\
&\ge &-\rho(\g(t),\g(s)) + \rho(\g(t),\g(0)) 
- \rho(\g(s),\g(0)) > -3\e,
\end{eqnarray*}
by part c) of Lemma \ref{lem4.4}.

Now let $m = \overline{\lim} \var_y(\g(t))$. 
Since $\var_y(x) \le \rho(y, \g(0))$ for all $x \in X$, we must
have $m \le \rho(y, \g(0))$. Now there is an $s_0 \ge N$ such
that $\var_y(\g(s_0)) \ge m - \e$. Set $M = s_0$. Then for
$t \ge M$ we must have $m \ge \var_y(\g(t)) \ge m -4\e$ according
to the previous paragraph. 
It follows that $\lim \var_y(\g(t)) = m$.

We can now show that $\g$ is weakly-geodesic. Given $\e > 0$, choose
$N$ and $M \ge N$ as above. Then for $t \ge s \ge M$ the
first condition of Definition 4.3c is satisfied by Lemma 4.4a,
while for the second condition we have from above
\begin{eqnarray*}
|\rho(\g(t), y) &- &\rho(\g(s), y) - (t - s)|\\
&\le &|\rho(\g(t), y) - \rho(\g(t), \g(0)) - 
\rho(\g(s), y) + \rho(\g(s), \g(0))|\\
& + &|\rho(\g(t), \g(0)) - t|
+ |\rho(\g(s), \g(0)) - s|\\
&\le &|\var_y(\g(t)) - \var_y(\g(s))| + 2\e < 6\e.
\end{eqnarray*}

Finally, suppose that $\g$ is a weakly-geodesic ray. For any $y \in X$
we show that $\{\var_y(\g(t))\}$ is a Cauchy net. Let $\e$ and $N$
be as in Definition 4.3c. Then for $t, s \ge N$ we have 
\begin{eqnarray*}
|\var_y(\g(t)) &- &\var_y(\g(s))|\\
&= &|\rho(\g(t), \g(0)) - \rho(\g(t), y) - 
\rho(\g(s), \g(0)) + \rho(\g(s), y)|\\
&\le &|\rho(\g(s), y) - \rho(\g(t), y) -(s-t)|\\
&+ &|\rho(\g(t), \g(0)) - t| + |s- \rho(\g(s), \g(0))| < 3\e.
\end{eqnarray*}
\end{proof}

For the next theorem we will need:

\setcounter{proposition}{5}
\begin{proposition}
\label{prop4.6}
Let $(X,\rho)$ be a locally compact metric space.  
If the topology of $X$ has a
countable base, then so do the topologies of ${\bar X}^{\rho}$ and 
$\p_{\rho}X$.
\end{proposition}

\begin{proof}
If $(X,\rho)$ is a locally compact metric space 
whose topology has a
countable base, then $C_{\infty}(X)$ has a countable dense set. 
Also, $X$
has a countable dense set, and the corresponding $\var_y$'s can be 
used to construct
a countable dense subset of $H_{\rho}$.  Thus the $C^*$-algebra ${\mathcal
G}(X,\rho)$ will have a countable dense set, and so the underlying
spaces will have countable bases for their topologies.
\end{proof}

We recall that a metric is said to be {\em proper} if every closed
ball of finite radius is compact.

\setcounter{theorem}{6}
\begin{theorem}
\label{th4.7}
Let $(X,\rho)$ be a complete locally compact metric space, and let $\g$ be a
weakly-geodesic 
ray in $X$.  Then $\lim_{t \rightarrow \infty} 
f(\g(t))$ exists for
every $f \in {\mathcal G}(X,\rho)$, and defines an element of $\p_{\rho}X$.
Conversely, if $\rho$ is proper and if the topology of $(X, \rho)$ has
a countable base, then every point of $\p_{\rho}X$ is determined as
above by a weakly-geodesic ray.
\end{theorem}

\begin{proof}
It is clear that the limit exists for the constant functions. 
From the definition of a weakly geodesic ray
we see that $\g$ must leave any compact set.  
Thus the limit
exists and is $0$ for all $f \in C_{\infty}(X)$.  Choose $\g(0)$ as 
the base-point
in defining $\var_y$ for any $y \in X$. Then from Lemma \ref{lem4.5}
we know that $\lim \var_y(\g(t))$ exists for all $y \in X$. 

Let ${\tilde {\mathcal G}}(X,\rho)$ denote the subalgebra of $C_b(X)$ 
generated by
$C_{\infty}(X)$, the constant functions, and the $\var_y$'s, before taking the
norm-closure.  It is clear from the above that $\lim f(\g(t))$ exists 
for every $f
\in {\tilde {\mathcal G}}(X,\rho)$, and that $|\lim f(\g(t))| \le 
\|f\|_{\infty}$.
Thus the limit defines a homomorphism from ${\tilde {\mathcal G}}(X,\rho)$ to
${\mathbb C}$ which is norm-continuous, and so extends to all of ${\mathcal
G}(X,\rho)$ by continuity.   
It thus defines a point, say $u$, of ${\bar X}^{\rho}$.  But
because $\g$ leaves any compact subset of $X$, the point defined by 
the limit must
be in $\p_{\rho}X$.  It is easy to check now that $\lim f(\g(t))$ 
exists and equals $f(u)$ for all $f
\in {\mathcal G}(X,\rho)$.

Suppose now that the topology of $(X,\rho)$ has a countable base, and
that $\rho$ is proper.
Let $u \in \p_{\rho}X$. Then we can apply Proposition \ref{prop4.6}
to conclude that there is a sequence, $\{w_n\}$, in $X$ which converges
in ${\bar X}^{\rho}$ to $u$. Since $u \notin X$ and $\rho$ is proper,
the sequence $\{w_n\}$ must be unbounded. Thus we can find a subsequence,
which we denote again by $\{w_n\}$, such that if $n > m$ then
$\rho(w_n, w_0) > \rho(w_m, w_0)$. Let $T$ denote the set of
$\rho(w_n, w_0)$'s, and for any $t \in T$ with $t = \rho(w_n, w_0)$
set $\g(t) = w_n$. Then $\lim \g(t) = u$. We show that $\g$ is
weakly-geodesic. Notice that by construction $\rho(\g(t), \g(0)) = t$
for each $t \in T$, so that the first condition of Definition 4.3c
is satisfied. Let $y \in X$. Use $\g(0)$ as the base-point
for defining $\var_y$. Now $\var_y(\g(t))$ converges to $\var_y(u)$,
and so, given $\e > 0$, we can find an $N$ such that whenever $s, t \in T$ 
with $s, t \ge N$ then $|\var_y(t) - \var_y(s)| \le \e$. 
Then for such $s, t$ we have
\[|\rho(\g(t), y) - \rho(\g(s), y) -(t - s)| = 
|\var_y(\g(t)) - \var_y(\g(s))| \le \e.
\]
\end{proof}

In view of the history of these ideas (see $1.2$ of \cite{G3}), we make:

\setcounter{definition}{7}
\begin{definition}
\label{def4.8}
A point of $\p_{\rho}X$ which is defined as above by an 
almost-geodesic ray $\g$
will be called a {\em Busemann point} of $\p_{\rho}X$, and we will 
denote the point
by $b_{\g}$.
\end{definition}

For any $(X, \rho)$ it is an interesting question as to whether
every point of $\p_{\rho}X$ is a 
Busemann point.
This is known to be the case for CAT(0) spaces (corollary II.8.20 of 
\cite{BrH}).
But in the next section we will need to deal with metric spaces which are not
CAT(0).  We will also see there by example that two metrics $\rho_1$ 
and $\rho_2$ on
$X$ which are Lipschitz equivalent, in the sense that there are 
positive constants
$k$, $K$ such that
\[
k\rho_1 \le \rho_2 \le K\rho_1,
\]
can give metric boundaries for $X$ which are not homeomorphic.

Here is an example of a complete locally compact non-compact metric 
space $X$ which
has no geodesic rays, but for which every point of $\p_{\rho}X$ is a Busemann
point.  Let $X$ be the subset $X = \{(n,1/n): n \ge 1\}$ of ${\mathbb 
R}^2$, with the
restriction to it of the Euclidean metric on ${\mathbb R}^2$.  This 
suggests the
usefulness of almost-geodesic rays.
Just before Proposition 5.4 we will give an example
of a proper metric on $\mathbb Z$ for which there are no almost-geodesic
rays, so no Busemann points (but there are sufficiently many 
weakly-geodesic rays).

We will later need:

\setcounter{proposition}{8}
\begin{proposition}
\label{prop4.9}
Let $z_0 \in X$ and let $\g$ and $\g'$ be almost-geodesic rays from 
$z_0$ (i.e.,
$\g(0) = z_0  = \g'(0)$).  If for any positive integer $N$ and any $\e > 
0$ we can find
$s$ and $t$ in the domains of $\g$ and $\g'$ respectively such that 
$s,t \ge N$ and
$\rho(\g(s),\g'(t)) < \e$, then $b_{\g} = b_{\g'}$.
\end{proposition}

\begin{proof}
Each $\var_y$ has Lipschitz constant $\le 2$, so
\[
|\var_y(\g(s)) - \var_y(\g'(t))| \le 2\rho(\g(s),\g'(t)).
\]
The desired result follows quickly from this.
\end{proof}

We now briefly consider isometries.  Suppose that $\a$ is an isometry of
$(X,\rho)$ onto itself.  
Then for $y,z \in X$ we have $\var_{y,z} \circ \a^{-1} =
\var_{\a(y),\a(z)}$.  Thus $H_{\rho}$ is carried onto itself by $\a$. 
Clearly so
are $C_{\infty}(X)$ and the constant functions, and so $\a$ gives an 
automorphism of
the algebra ${\mathcal G}(X,\rho)$.  It follows that $\a$ gives a 
homeomorphism of
${\bar X}^{\rho}$ onto itself which extends $\a$ on $X$.  This 
homeomorphism carries
$\p_{\rho}X$ onto itself.  Thus:

\setcounter{proposition}{9}
\begin{proposition}
\label{prop4.10}
Every isometry of a complete locally compact metric space $(X,\rho)$ extends
uniquely to a homeomorphism of ${\bar X}^{\rho}$ onto itself which carries
$\p_{\rho}X$ onto itself.
\end{proposition}

Later we will need to consider (cartesian) products of metric spaces. 
There are
many ways to define a metric on a product.  One of these ways meshes especially
simply with the construction of the metric compactification.  If 
$(X,\rho_X)$ and
$(Y,\rho_Y)$ are metric spaces, we define $\rho$ on $X \x Y$ by
\[
\rho((x_1,y_1),(x_2,y_2)) = \rho_X(x_1,x_2) + \rho_Y(y_1,y_2).
\]
We will call $\rho$ the ``sum of metrics''.

\setcounter{proposition}{10}
\begin{proposition}
\label{prop4.11}
Let $(X,\rho_X)$ and $(Y,\rho_Y)$ be locally compact metric spaces, 
and let $\rho$
be the sum of metrics on $X \x Y$.  Then
\[
(X \x Y)^{-\rho} = ({\bar X}^{\rho_X}) \x ({\bar Y}^{\rho_Y}).
\]
\end{proposition}

\begin{proof}
We need to show that the evident map from $X \x Y$ to 
$({\bar X}^{\rho_X}) \x ({\bar Y}^{\rho_Y})$ extends to a
homeomorphism from $(X \x Y)^{-\rho}$. For this it suffices to show that
the restriction map from $C(({\bar X}^{\rho_X}) \x ({\bar Y}^{\rho_Y}))$
to $C_b(X \x Y)$ maps into $C((X \x Y)^{-\rho})$ and is onto.
Let $x_0,y_0$ be base-points in $X$ and $Y$ respectively, and use 
$(x_0,y_0)$ as a
base-point for $X \x Y$.  Then for $(u,v) \in X \x Y$ we have
\begin{eqnarray*}
\var_{(u,v)}(x,y) &= &\rho((x,y),(x_0,y_0)) - \rho((x,y),(u,v)) \\
&= &\rho_X(x,x_0) - \rho_X(x,u) + \rho_Y(y,y_0) - \rho_Y(y,v) \\
&= &\var_u(x) + \var_v(y).
\end{eqnarray*}
In particular, $\var_{(u,y_0)} = \var_u \otimes 1_Y$ and 
$\var_{(x_0,v)} = 1_X \otimes \var_v$. Thus the restrictions of
$\var_u \otimes 1_Y$ and $1_X \otimes \var_v$ are in $C((X \x Y)^{-\rho})$.
The same is true for any $f \otimes 1_Y$ and $1_X\otimes g$
where $f \in C_c(X)$ and $g \in C_c(Y)$, or for constant functions. Thus
the range of the restriction map is in $C((X \x Y)^{-\rho})$. But from
the calculation above we also see that any $\var_{(u,v)}$ is in the range of
the restriction map, and from this it is easily seen that the restriction
map is onto $C((X \x Y)^{-\rho})$.
\end{proof}

%%%%%%%%%%%%%%%%%%%%%%%%%%%%%%%%%%%%%%%%%%%%%%%%%%%%
\setcounter{section}{4}
\section{The case of $G = {\mathbb Z}$}
\label{sec5}
%%%%%%%%%%%%%%%%%%%%%%%%%%%%%%%%%%%%%%%%%%%%%%%%%%%%%

In this section we will see how the constructions of the previous 
sections can be
used to deal with Questions \ref{quest2.6} and \ref{quest2.7} when $G 
= {\mathbb
Z}$.  This case already reveals some phenomena which we will have to 
deal with later
for the case $G = {\mathbb Z}^d$.

\setcounter{example}{0}
\begin{example}
\label{exam5.1}
{\em We examine first the case in which $\ell$ is the standard length 
function on $G =
{\mathbb Z}$ defined by $\ell(n) = |n|$, so that $\rho(m,n) = |m-n|$.  
Note that $\ell$ is the word-length function for the
generating set $S = \{\pm 1\}$.  We determine $\p_{\ell}G$.  For any 
$k \in {\mathbb
Z}$ we have
\[
\var_k(n) = |n| - |n-k|.
\]
In particular,
\[
\var_k(n) = \left\{ \begin{array}{rl}
k &\mbox{for $n \ge 0$ and $n \ge k$} \\
-k &\mbox{for $n \le 0$ and $n \le k$.}
\end{array} \right.
\]
 From this it is clear that ${\bar {\mathbb Z}}^{\ell}$ is just 
${\mathbb Z}$ with
the points $\{\pm \infty\}$ adjoined in the traditional way.  The 
action $\a$ of
${\mathbb Z}$ on ${\bar {\mathbb Z}}^{\ell}$ is by translation 
leaving the points at
infinity fixed.  Thus $\p_{\ell}Z = \{\pm \infty\}$ 
with the trivial action $\a$ of
${\mathbb Z}$.

Now let $f \in C_c({\mathbb Z})$ be given.  Since ${\mathbb Z}$ is 
amenable, we know
that $[M_{\ell},\pi(f)]$ is in $C({\bar {\mathbb Z}}^{\ell}) \x_{\a} 
{\mathbb Z}$,
and that this crossed product is faithfully represented on 
$\bell^2({\mathbb Z})$,
as discussed in Section~\ref{sec3}.  We can factor by ${\mathcal K} =
C_{\infty}({\mathbb Z}) \x_{\a} Z$, and so look at the image of 
$[M_{\ell},\pi_f]$
in the cosphere algebra $S^*A$, which by the discussion of 
Section~\ref{sec3} is
exactly $C(\p_{\ell}{\mathbb Z}) \x_{\a} {\mathbb Z}$.  This latter 
is isomorphic to
two copies of $C^*({\mathbb Z})$.  The image of $\Sigma 
f(y)M_{\var_y}\pi_y$ in the
copy at $+\infty$ will be $\{k \mapsto kf(k)\}$, 
while the image in the copy at
$-\infty$ will be $\{k \rightarrow -kf(k)\}$.  Let us take here the 
convention that
the Fourier series for any $g \in C_c({\mathbb Z})$ is given by ${\hat g}(t) =
\Sigma g(k)e^{ikt}$, so that ${\hat g}'(t) = i \Sigma kg(k)e^{ikt}$. 
Then we see
from just above that
\[
L(f) = \|\Sigma 
f(y)M_{\var_y}\pi_y\| \ge \|{\hat f}'\|_{\infty}.
\]
But $\|{\hat f}'\|_{\infty}$ agrees with the standard Lip-norm on 
$C^*({\mathbb Z})
= C({\mathbb T})$ which gives the circle a circumference of $2\pi$.  From the
comparison lemma
$1.10$ of \cite{R4} it follows that $L$ is a Lip-norm, and that it 
gives ${\mathbb
T}$ (and so the state space $S(C^*(Z))$) radius no larger than $\pi$. }
\end{example}

\setcounter{example}{1}
\begin{example}
\label{exam5.2}
{\em Again we take $G = {\mathbb Z}$, but now we take the word-length 
function $\ell$
corresponding to the generating set $\{\pm 1,\pm 2\}$.  
Then $\ell$ is given by
\[
\ell(n) = [|n|/2], 
\]
where $[\cdot]$ denotes ``least integer not less than''.  Thus for any $k \in
{\mathbb Z}$
\[
\var_k(n) = [|n|/2] - [|n-k|/2].
\]
 From this one finds that if $k$ is even then
\[
\var_k(n) = \left\{ \begin{array}{rl}
k/2 &\mbox{for $n \ge 0$ and $n \ge k$} \\
-k/2 &\mbox{for $n \le 0$ and $n \le k$,}
\end{array} \right.
\]
whereas if $k$ is odd then
\[
\var_k(n) = \left\{ \begin{array}{ll}
\left. \begin{array}{l}
(k-1)/2 \mbox{ for $n$ even} \\
(k+1)/2 \mbox{ for $n$ odd}
\end{array} \right\} & \mbox{ for } n \ge 0 \mbox{ and } n \ge k \\
\left. \begin{array}{l}
-(k+1)/2 \mbox{ for $n$ even} \\
-(k-1)/2 \mbox{ for $n$ odd}
\end{array} \right\} & \mbox{ for }n \le 0 \mbox{ and } n \le k.
\end{array} \right.
\]
 From this it is easily seen that $\p_{\ell}{\mathbb Z}$ will consist 
of $4$ points,
two at $+\infty$ and two at $-\infty$, which we can label ``even'' 
and ``odd''.  The
action of ${\mathbb Z}$ on $\p_{\ell}{\mathbb Z}$ will at each end be that of
${\mathbb Z}$ on ${\mathbb Z}_2 = {\mathbb Z}/2{\mathbb Z}$.  In 
particular, the
boundary contains no fixed-points for this action. }
\end{example}

We learn several things from comparing this example with the one just before.
First, two word-length metrics on a given group can give metric 
boundaries which are
not homeomorphic.  But it is well-known (e.g., proposition $8.3.18$ 
of \cite{BBI})
and easily seen that if $G$ is a finitely-generated group and if $\ell_1$ and
$\ell_2$ are the word-length functions for two finite generating 
sets, then the
corresponding left-invariant metrics are (Lipschitz) equivalent in 
the sense defined
in the previous section.  Thus we see that equivalent metrics which 
give (the same)
locally compact topologies (even discrete) and for which the set is 
complete, can
give metric boundaries which are not homeomorphic.

Next, ${\mathbb Z}$ is an example of a hyperbolic group \cite{GhH}, 
and so for the
metric from either of these generating sets it is a hyperbolic metric 
space.  But
the Gromov boundary of a hyperbolic space is independent of the 
metrics as long as
the metrics are equivalent, or at least coarsely equivalent.  The 
Gromov boundary
for ${\mathbb Z}$ is just $\{\pm \infty\}$.  One way of viewing what 
is happening is
that for the metric of the Example 5.2 the maps $m \mapsto 2m$ 
and $m \mapsto
2m+1$ are geodesic rays which determine Busemann points in the 
boundary which
are our two points at $+\infty$.  But for the Gromov boundary any two 
geodesic rays
which stay a bounded distance from each other define the same point 
at infinity.
In particular, our present example shows that for a given hyperbolic 
metric space
the metric boundary and the Gromov boundary can fail to be homeomorphic.

For our next observation, let $(X,\rho)$ be a proper metric space with 
base-point $z_0$,
and let $T \subset {\mathbb R}^+$ be a fixed domain for geodesic 
rays, so that $0
\in T$ and $T$ is unbounded.  On the set of
geodesic rays from $z_0$ whose domain is $T$ we put the 
topology of pointwise convergence (which, because geodesic rays
are Lipschitz maps of Lipschitz constant $1$, is equivalent to the
topology of
uniform convergence on bounded subsets of $T$).  This is done in 
various places in
the literature.  Because $\rho$ is proper, it is easy to 
see that the set
of all such geodesic rays is compact. 
For groups $G$ with a word-length $\ell$ (or for graphs in
general) it is natural to take $T = {\mathbb Z}^+$. 
It is reasonable to wonder then whether 
$\p_{\ell}G$ is
the quotient of this compact set of geodesics, with the quotient 
topology.  If it
were, then for each $y \in G$ the function which assigns to each such
geodesic ray $\g$ from $e$ the number $\lim\var_y(\g(t))$
should be a continuous 
function on this
compact set.  But this already fails for Example 5.2.  
For each $k \ge 1$ let $\g^k$ be the 
geodesic ray
from $0$ defined by
\[
\g^k(n) = \left\{ \begin{array}{rl}
2n &\mbox{if $n \le k$} \\
2n-1 &\mbox{if $n \ge k+1$.}
\end{array} \right.
\]
Then $\g^k$ converges pointwise to the geodesic ray defined by 
$\g^{\infty}(n) =
2n$ for all $n$.  But it is easy to see that $b_{\g^{\infty}}$ is the 
even point at
$+\infty$ while $b_{\g^k}$ is the odd point at $+\infty$ for all $k$.  We also
remark that in our present example there is no geodesic line which 
joins the two
points at $+\infty$ (so this example fails to have the property of 
``visibility''
\cite{GhH}).

Our Example 5.2 also shows that the metric compactification is 
not in general
well-related to the Higson compactification, 
as defined in $5.4$ of \cite{Roe}.
For that definition let $(X,\rho)$ be a proper metric space.  For any 
$r > 0$ we
define the variational function, $V_rf$, of any function $f$ by
\[
(V_rf)(x) = \sup\{|f(x)-f(y)|: \rho(y,x) \le r\}.
\]
The Higson compactification is the maximal ideal space of the unital 
commutative
$C^*$-algebra of all bounded continuous functions on $X$ such that 
for each $r > 0$
the function $V_rf$ vanishes at infinity.
For Example 5.2 let us consider $V_2\var_1$.  Easy 
calculation shows that
for any $n \ge 1$ we have $\var_1(2k) = 0$ while $\var_1(2k+1) = 1$.  But
$\rho(2k,2k+1) = \ell(1) = 1$ for all $k$.  Thus $(V_2\var_1)(k) \ge 
1$ for all
$k$.  Consequently $\var_1$ does not extend to the Higson 
compactification.  More
generally, if a complete locally compact metric space $(X,\rho)$ has 
geodesic rays
which determine distinct Busemann points of $\p_{\rho}X$ and yet stay a finite
distance from each other, then ${\bar X}^{\rho}$ is not a quotient of 
the Higson
compactification.  Indeed, since the $\var_y$'s separate the points 
of $\p_{\rho}X$,
there will be some $y$ such that its $\var_y$ separates the two 
Busemann points, and
$V_r\var_y$ will not vanish at infinity if $r$ is larger than the 
distance between
the two rays.

The situation becomes yet more interesting when we consider 
generating sets such as
$\{\pm 3,\pm 8\}$.  But the proof given above that we obtain a 
Lip-norm when we use
the generating set $\{\pm 1\}$ extends without too much difficulty to 
the case of
arbitrary finite generating sets for ${\mathbb Z}$.  We do not 
include this proof
here since in Section~\ref{sec9} we will treat the general case of 
${\mathbb Z}^d$
by similar techniques, though the details are certainly more complicated.

However we will discuss here another 
approach for the case of $G = {\mathbb Z}$
which uses a classical argument which was pointed out to me by 
Michael Christ.  (I thank him for his guidence in this matter). This
second approach seems less likely to generalize to more complicated 
groups, but it
gives a stronger result for ${\mathbb Z}$.  For any $\b$ with $0 < \b 
\le 1$ and any
metric $\rho$ on a set, $\rho^{\b}$ will again be a metric, because 
$t \rightarrow
|t|^{\b}$ is a length function on ${\mathbb R}$.  In 
particular, if we set
$\ell_{\b}(n) = |n|^{\b}$ then $\ell_{\b}$ is a length function on 
${\mathbb Z}$.

\setcounter{theorem}{2}
\begin{theorem}
\label{th5.3}
Let $\o$ be a translation-bounded function on ${\mathbb Z}$ such that 
$\o(0) = 0$.
If $\ell_{\b}/\o$ is a bounded function (ignoring $n = 0$) for some 
$\b$ with $1/2
< \b \le 1$, then $L_{\o}$ is a Lip-norm on $C^*({\mathbb Z}) = 
C({\mathbb T})$.
\end{theorem}

\begin{proof}
For any group $G$ and any $\o$ we have
\[
[M_{\o},\pi_f]\d_e = \Sigma f(y)\var_y(y)\d_y = \Sigma \o(y)f(y)\d_y,
\]
where $\{\d_y\}$ here denotes the standard basis for $\bell^2(G)$.  Thus
\[
\|\o f\|_2 \le \|[M_{\o},\pi_f]\| = L_{\o}(f).
\]
What is special about ${\mathbb Z}$ is that $\|\o f\|_2$ can control 
the norms we
need.  (This is related to our discussion of ``rapid decay''
in Section 1.)  For this we need that
$\ell_{\b}^{-1} \in \bell^2({\mathbb Z})$, which happens exactly for 
$\b > 1/2$.
(Here and below we ignore $n=0$ or set $\ell_{\b}^{-1}(0) = 0$.)  Let
$f \in C_c({\mathbb Z})$, with ${\hat f}$ its Fourier transform on 
${\mathbb T} = {\mathbb R}/{\mathbb Z}$, viewed as a periodic
function on $\mathbb R$.  For  $s, t \in {\mathbb R}$ with $s < t$
and $|t-s| < 1$ let 
$\chi_{[s,t]}$ denote the
characteristic function of the interval $[s,t]$, extended 
by periodicity.  Then
%\begin{eqnarray*}
\[
|{\hat f}(s)-{\hat f}(t)| = \left| \int_s^t {\hat f}'(r)dr\right| 
= |\<{\hat f}',\chi_{[s,t]}\>| 
= |\<({\hat f}')^{\vee},(\chi_{[s,t]})^{\vee}\>|.
\]
%\end{eqnarray*}
But $(\chi_{[s,t]})^{\vee}(n) = (1/i2\pi n)(e(nt) - e(ns))$ if we
set $e(r) = e^{2\pi inr}$, while $({\hat f}')^{\vee}(n) = -2\pi 
inf(n)$.  Thus if
we set $g_{s,t}(n) = (e(nt) - e(ns))$, the above becomes
$|\<f,g_{s,t}\>|$ as a pairing between functions in $\bell^1({\mathbb Z})$ and
$\bell^{\infty}({\mathbb Z})$.  But (with $\o^{-1}(0) = 0$) we can 
rewrite this as
\[
|\<\o f,\o^{-1}g_{s,t}\>| \le \|\o f\|_2 \|\o^{-1}g_{s,t}\|_2,
\]
and notice that
\begin{eqnarray*}
\|\o^{-1}g_{s,t}\|_2 &= &\|(\ell_{\b}/\o)\ell_{\b}^{-1}g_{s,t}\|_2 \\
&\le &\|\ell_{\b}/\o\|_{\infty} \|\ell_{\b}^{-1} g_{s,t}\|_2 < \infty,
\end{eqnarray*}
since $\ell_{\b}^{-1} \in \bell^2({\mathbb Z})$.  Set $m(s,t) = 
\|\ell_{\b}^{-1}
g_{s,t}\|_2$.  Then putting the above together, we obtain
\[
|{\hat f}(t) - {\hat f}(s)| \le m(s,t) \|\ell_{\b}/\o\|_{\infty} 
\|L_{\o}(f)\|.
\]
  A simple estimate using 
the fact that
$\ell_{\b}^{-1} \in \bell^2({\mathbb Z})$ shows that for each $\e > 
0$ there is a
$\d > 0$ such that if $|t-s| < \d$ then $m(s,t) < \e$.  From this we 
see that the
set of $\hat f$'s for which $L_{\o}(f) \le 1$ and $f(0) = 0$ forms a 
bounded subset of
$C({\mathbb T})$ which is equicontinuous, so totally bounded by the 
Arzela--Ascoli
theorem.  From this it is clear that $L_{\o}$ gives finite radius 
and, by theorem
$1.9$ of \cite{R4}, that it is a Lip-norm.
\end{proof}

I suspect that when $\b < 1/2$ then $L_{\ell_{\b}}$ fails to be a 
Lip-norm, but I
have not found a proof of this.

Notice that Theorem 5.3 applies if $|\o(n)| \ge 1$ for $n \neq 0$
and if there are positive constants $c$ and $K$ such that
$|\o - c\ell^\b| \le K$, for then $|\ell/\o| \le (K+1)c$. This is
the situation which occurs for the various word-length functions
on $\mathbb Z$ (for $\b = 1$).

It is interesting to see what the metric compactification of 
${\mathbb Z}$ is when
$\b < 1$.  For any $p \in {\mathbb Z}$ we have
\[
\var_p(n) = |n|^{\b} - |n-p|^{\b} = \int_{|n-p|}^{|n|}\b t^{\b-1}\ dt.
\]
Since $t^{\b-1} \rightarrow 0$ at $+\infty$ because $\b < 1$, it follows that 
$\var_p(n)
\rightarrow 0$ as $n \rightarrow \pm \infty$.  Thus $\var_p \in
C_{\infty}({\mathbb Z})$, and so the metric compactification is just 
the one-point
compactification of ${\mathbb Z}$.  Note also that 
$[M_{\ell_{\b}},\pi_f]$ is a
compact operator for each $f \in C_c({\mathbb Z})$.  Thus the 
cosphere algebra for
$(C^*({\mathbb Z}),\bell^2({\mathbb Z}),M_{\ell_{\b}})$ is $C^*(\mathbb Z)$,
and the image of $[M_{\ell_\b}, \pi_f]$ in it is $0$.
We also remark that it is easily verified that if we set
$\g(n^{\b}) = n$, then $\g$ is a weakly-geodesic ray, but that
there are no almost-geodesic rays in $\mathbb Z$ for this metric,
since by parts a) and b) of Lemma 4.4 if $\g$ were such a ray we would
have, for any fixed big $r$, that $|\g(t)|^\b - |\g(t-r)|^b$ would be 
approximately $r$ as $t \rightarrow \infty$, contradicting our observation
above that it must go to 0.

We conclude this section with the following observation, which 
applies to our more
general case of ${\mathbb Z}^d$.

\setcounter{proposition}{3}
\begin{proposition}
\label{prop5.4}
Let $\o$ be a translation-bounded function on a countable discrete 
Abelian group
$G$, let $L_{\o}$ on $C_c(G)$ be defined as earlier by $L_{\o}(f) = 
\|[M_{\o},\pi_f]\|$, and
let $\rho_{\o}$ be the corresponding metric on ${\hat G}$ (which may 
not give the
usual topology of ${\hat G}$).  Then $\rho_{\o}$ is invariant under 
translation on
${\hat G}$.
\end{proposition}

\begin{proof}
Let us denote the pairing between $G$ and ${\hat G}$ by $\<m,t\>$. 
Then translation
on ${\hat G}$ corresponds to the dual action, $\b$, of ${\hat G}$ on 
$C^*(G)$ given
on $C_c(G)$ by $(\b_t(f))(m) = \<m,t\>f(m)$.  This is unitarily implemented in
$\bell^2(G)$ by $M_t$, where $(M_t\xi)(m) = \<m,t\>\xi(m)$.  Then
\begin{eqnarray*}
[M_{\o},\b_t(\pi_f)] &= &[M_{\o},M_t\pi_fM_t^*] \\
&= &M_t[M_{\o},\pi_f]M_t^*,
\end{eqnarray*}
so that $L_{\o}(\b_t(f)) = L_{\o}(f)$.  (In other words, $\b$ is an action by
isometries as defined in \cite{R6}.)
\end{proof}

 From Theorem 5.3 one begins to see that ${\mathbb T}^d$ 
has a bewildering
variety of translation invariant metrics which give its topology. 
For example, if
$\rho$ is such a metric then so is $\rho^r$ for any $r$ with $0 < r < 
1$, as is any
convex function of $\rho$.  The sum of two metrics and the supremum 
of two metrics
are again metrics.  More generally, the ``$\bell^p$-sum'' of two metrics is a
metric.  These operations all preserve translation invariance.  
For the case of
${\mathbb T}$, any strictly increasing continuous function $\ell$ on 
$[0,1/2]$ such
that $\ell(0) = 0$ and $\ell(s+t) \le \ell(s) +  \ell(t)$ if $s+t 
\le 1/2$ gives in
an evident way a continuous length function on ${\mathbb T} = 
{\mathbb R}/{\mathbb
Z}$, and all continuous length functions on ${\mathbb T}$ arise in 
this  way.  It
would be interesting to determine which generating sets for ${\mathbb 
Z}$ determine
which length functions on ${\mathbb T}$, but I have not investigated 
this question.

%%%%%%%%%%%%%%%%%%%%%%%%%%%%%%%%%%%%%%%%%%%%
\setcounter{section}{5}
\section{The metric compactification for norms on ${\mathbb R}^d$}
\label{sec6}
%%%%%%%%%%%%%%%%%%%%%%%%%%%%%%%%%%%%%%%%%%%

One of our eventual aims is to show that when $\ell$ is a length function on
${\mathbb Z}^d$ which is the restriction to ${\mathbb Z}^d$ of a norm 
on ${\mathbb
R}^d$, then $L_{\ell}$ is a Lip-norm.  In preparation for this we 
examine here the
metric compactification of ${\mathbb R}^d$ for any given norm.  We begin by
considering the usual $\bell^1$-norm, both because it is simple to treat and
displays some interesting phenomena, and also because its restriction 
to ${\mathbb
Z}^d$ gives the word-length function for the standard generating set. 
Following up
on Example \ref{exam5.1}, we set ${\bar {\mathbb R}} = {\mathbb R} \cup \{\pm
\infty\}$ in the usual way, with the action of ${\mathbb R}$ fixing 
the points $\pm
\infty$.

\setcounter{proposition}{0}
\begin{proposition}
\label{prop6.1}
The metric compactification of $({\mathbb R}^d,\|\cdot\|_1)$ is just $({\bar
{\mathbb R}})^d$ with its product action of ${\mathbb R}^d$.  Thus the metric
boundary is the set of $({\tilde x}_j) \in ({\bar {\mathbb R}})^d$ such that at
least one entry is $+\infty$ or $-\infty$.
\end{proposition}

\begin{proof}
The metric from $\|\cdot\|_1$ on ${\mathbb R}^d$ is easily seen to be 
the sum of the
metrics on ${\mathbb R}$ in the sense used in Proposition 
\ref{prop4.11}.  Thus we
just need to apply that proposition a number of times.
\end{proof}

We note that now there are orbits in the boundary which are not 
finite, but there
are also fixed points (only a finite number of them).

We now investigate what happens for other norms on ${\mathbb R}^d$.  It is
notationally convenient for us just to consider a finite-dimensional 
vector space
$V$ with some given norm $\|\cdot\|$.  We will denote the corresponding metric
boundary simply by $\p_{\ell}V$, where $\ell(x) = \|x\|$ for all $x \in V$.

For any $v \in V$ with $\|v\| = 1$ it is evident that the function 
$\g(t) = tv$ for
$t \in T = [0,\infty)$ is a geodesic ray, and so from our earlier 
discussion it will
determine a Busemann point, $b_v$, in $\p_{\ell}V$.  We now convert 
to this picture
some of the known elementary facts about 
tangent functionals of convex sets, as
explained for example in section V.9 of \cite{DfS}.  There is at 
least one linear
functional, say $\s$, on $V$ such that $\|\s\| = 1 = \s(v)$.  We call 
such a $\s$ a
``support functional'' at $v$.  Then for any $y \in V$ we have
%\begin{eqnarray*}
\[
\var_y(\g(t)) = \|tv\| - \|tv-y\| 
\le t - \s(tv-y) = \s(y).
\]
%\end{eqnarray*}
In particular, $\var_{-y}(\g(t)) \le -\s(y)$.  On letting $t$ go to 
$+\infty$ we
find that
\[
-\var_{-y}(b_v) \ge \s(y) \ge \var_y(b_v).
\]
But theorem~$5$ of section V.9 of \cite{DfS} (which uses the 
Hahn--Banach theorem)
tells us that for any real number $r$ such that $-\var_{-y}(b_v) \ge r \ge
\var_y(b_v)$ there is a support functional $\s$ at $v$ such that 
$\s(y) = r$.  To
see that theorem V.9.5 really applies here, we note that if we set $s 
= t^{-1}$ then
\[
\|tv\| - \|tv-y\| = (\|v\| - \|v-sy\|)/s,
\]
and that $s \rightarrow +0$ as $t \rightarrow +\infty$.  From this 
viewpoint we are
thus looking at the negative of the tangent functional to the unit 
ball at $v$ in
the direction of $-y$, which fits the setting of theorem V.9.5.

The point $v$ is called a {\em smooth} point of the unit sphere if 
there is only one
support functional $\s$ at $v$.  We denote this unique $\s$ by 
$\s_v$.  Then the
above considerations tell us that if $v$ is smooth then $\var_y(b_v) =
-\var_{-y}(b_v)$.  On combining this with the inequalities found 
above, we obtain:

\setcounter{proposition}{1}
\begin{proposition}
\label{prop6.2}
Let $v$ be a smooth point of the unit sphere of $V$.  Then
\[
\var_y(b_v) = \s_v(y)
\]
for all $y \in V$.
\end{proposition}

For us the following proposition will be of considerable importance. 
We consider
the action of $V$ on itself by translation, and the corresponding action on
$\p_{\ell}V$.

\setcounter{proposition}{2}
\begin{proposition}
\label{prop6.3}
Let $v$ be a smooth point of the unit sphere of $V$.  
Then $b_v$ is a fixed point
under the action of $V$ on $\p_{\ell}V$.
\end{proposition}

\begin{proof}
We use the $1$-cocycle relation $2.2$ and Proposition \ref{prop6.2}
to calculate that for any $x,y 
\in V$ we have
\begin{eqnarray*}
(\a_x\var_y)(b_v) &= &\var_{x+y}(b_v) - \var_x(b_v) \\
&= &\s_v(x+y) - \s_v(x) = \s_v(y) = \var_y(b_v).
\end{eqnarray*}
\end{proof}

Finally, we note that theorem~$8$ of section V.9 of \cite{DfS} says 
that, for any
norm, the set of smooth points of the unit sphere is dense in the 
unit sphere.  This
does not imply that there are infinitely many fixed points in 
$\p_{\ell}V$, as the
next example shows.  But we will see later that it does show that 
there are enough
for our purposes.

\setcounter{example}{3}
\begin{example}
\label{exam6.4}
{\em We examine the case of ${\mathbb R}^2$ with $\|\cdot\|_1$, whose metric
compactification is described by Proposition \ref{prop6.1}.
Let us see how our
considerations concerning geodesics fit this example.  We identify 
the dual space
$V'$ in the usual way with ${\mathbb R}^2$ with the norm
\[
\|(r,s)\|_{\infty} = \max\{|r|,|s|\}.
\]

All but $4$ points of the unit sphere of $V$ are smooth.  
However, for any $v =
(a,b)$ with $0 < a$, $0 < b$ and $a+b=1$ we see that $\s_v = (1,1) 
\in V'$.  Thus
all these different $v$'s determine the same Busemann point of 
$\p_{\ell}V$.  This
accords with Proposition \ref{prop6.2} and the fact that $\p_{\ell}V$ 
has only $4$
fixed-points for the action of ${\mathbb R}^2$.

If instead we let $v$ be the non-smooth point $(1,0)$ and let $\g$ be the
corresponding geodesic ray, then for any $y = (p,q) \in {\mathbb R}^2$ we have
\begin{eqnarray*}
\var_y(\g(t)) &= &\|\g(t)\| - \|\g(t) - (p,q)\| \\
&= &|t| - |t-p| - |q|.
\end{eqnarray*}
The limit as $t \rightarrow +\infty$ is clearly $p-|q|$, so that
\[
\var_{(p,q)}(b_v) = p - |q| = \var_p(+\infty) + \var_q(0),
\]
where $\var_p$ and $\var_q$ are for ${\mathbb R}$.  Thus $b_v = 
(+\infty,0)$ in the
description of $\p_{\ell}V$ given by Proposition \ref{prop6.1}. 
Clearly $b_v$ is
not given by an element of $V'$.  It is easily seen that this $b_v$ 
is not invariant
under translation.

We see in this way that the {\em linear} geodesic rays from $0$, 
corresponding to
the points of the unit sphere, determine only $8$ Busemann points of 
$\p_{\ell}V$.
But we can show that every point of $\p_{\ell}V$ is determined by at least one
(possibly non-linear) geodesic ray from $0$.  For example, if we consider
$(+\infty,s) \in \p_{\ell}V$ for some fixed $s \in {\mathbb R}$, we 
can pick any
$t_0 \ge 0$ and let $\g$ consist of the unit-speed straight-line path 
from $(0,0)$
to $(t_0,0)$, followed by that from $(t_0,0)$ to $(t_0,s)$, followed 
by the linear
ray from $(t_0,s)$ in the direction $(1,0)$.  (We deal here with the 
``Manhattan
metric''.)  It is easy to check that $\g$ is a geodesic ray whose 
Busemann point
corresponds to $(+\infty,s)$.  We see in this way that every point of 
$\p_{\ell}V$
is a Busemann point.  It is also easy to see that for each of the $4$ 
points $(\pm
\infty,0)$ and $(0,\pm \infty)$ of $\p_{\ell}V$ there is only one 
geodesic ray to
them from $0$, but that for every other point of $\p_{\ell}V$ there 
are uncountably
many geodesic rays to it from $0$.}
\end{example}

\setcounter{question}{4}
\begin{question}
\label{quest6.5}
Is it true that, for every finite-dimensional vector space and every 
norm on it,
every point of $\p_{\ell}V$ is a Busemann point?
\end{question}

One says that $(V,\|\cdot\|)$ is {\em smooth} if every point of the 
unit sphere,
$S$, of $V$ is a smooth point.  Let $S'$ denote the unit sphere of 
$V'$.  Then our
earlier mapping $v \mapsto \s_v$ is defined on all of $S$. 
Furthermore it is onto
$S'$, because $V$, being finite dimensional, is reflexive.  This 
mapping $\s$ can
also be seen to be continuous.  This is essentially the fact that, as 
remarked at
the bottom of page~$60$ of \cite{LnT}, a compactness argument shows 
that smoothness
implies uniform smoothness.  However, if $S$ has ``flat spots'' then 
$\s$ will not
be injective.  It is not difficult to show that for $(V,\|\cdot\|)$ smooth,
$\p_{\ell}V$ can be naturally identified with $S'$, glued at $\infty$ 
using $\s$.
In this case each point of $\p_{\ell}V$ will be fixed by the action of $V$.

\setcounter{question}{5}
\begin{question}
\label{quest6.6}
For a general $(V,\|\cdot\|)$ is there an attractive description of 
$\p_{\ell}V$ and
of the action of $V$ on it?
\end{question}

We have seen in Example \ref{exam6.4} that the number of support 
functionals $\s_v$
coming from smooth points $v$ of the unit sphere can be finite.  The 
reason that
they nevertheless are adequate for our later 
purposes is given by the following
proposition (which must be already known):

\setcounter{proposition}{6}
\begin{proposition}
\label{prop6.7}
Let $\|\cdot\|$ be a norm on a finite-dimensional vector space $V$. 
Let $w \in V$,
and suppose that $|\s_v(w)| \le r$ for all smooth points $v$ of the unit 
sphere.  Then
$\|w\| \le r$.  Furthermore, the closed convex hull of $\{\s_v: v 
\mbox{ is a smooth
point}\}$ is the unit ball in the dual space $V'$ for the dual norm 
$\|\cdot\|'$.
\end{proposition}

\begin{proof}
Let $\|w\| = s$.  Because the smooth points are dense in the unit sphere by
theorem~$8$ of section V.9 of \cite{DfS}, 
for any $\e > 0$ we can find a smooth
point $v$ such that $\|w - sv\| < \e$.  
Then $|\s_v(w)-s| = |\s_v(w-sv)| < \e$.
Since $|\s_v(w)| \le r$ and $\e$ is arbitrary, 
it follows that $\|w\| = s \le r$.

Suppose now that $\tau \in V'$ and that $\tau \notin {\bar 
{co}}\{\s_v: v \mbox{
smooth}\}$.  Then by the Hahn--Banach theorem there is a $w \in V$ 
and an $r \in
{\mathbb R}$ such that $|\s_v(w)| \le r < \tau(w)$ for all smooth 
$v$.  But we have
just seen that then $\|w\| \le r$.  Thus $\|\tau\|' > 1$.
\end{proof}

%%%%%%%%%%%%%%%%%%%%%%%%%%%%%%%%%%%%%%%%%%%%%%%
\setcounter{section}{6}
\section{Restrictions of norms to ${\mathbb Z}^d$}
\label{sec7}
%%%%%%%%%%%%%%%%%%%%%%%%%%%%%%%%%%%%%%%%%%%%%%%%%%%%%%%

In this section we will examine what happens when norms on $V = 
{\mathbb R}^d$ are
restricted to ${\mathbb Z}^d$.  We begin with the case of the norm 
$\|\cdot\|_1$.
Following up on Example \ref{exam6.4} 
we set ${\bar {\mathbb Z}} = {\mathbb Z} \cup \{\pm
\infty\}$ in the usual way, with its action of ${\mathbb Z}$ leaving fixed the
points at infinity.  The proof of the following proposition is 
basically the same as
that of Proposition \ref{prop6.1}.

\setcounter{proposition}{0}
\begin{proposition}
\label{prop7.1}
For $\ell = \|\cdot\|_1$, the metric compactification of $({\mathbb 
Z}^d,\ell)$ is
$({\bar {\mathbb Z}})^d$ with its product action of ${\mathbb Z}^d$.  The 
metric boundary  is
the set of $({\tilde n}_j) \in ({\bar {\mathbb Z}})^d$ such that at 
least one entry
is $+\infty$ or $-\infty$.
\end{proposition}

Suppose now that $\ell = \|\cdot\|$ is any norm on $V = {\mathbb 
R}^d$, and that we
restrict it to ${\mathbb Z}^d$.  For any $y \in {\mathbb Z}^d$ the 
function $\var_y$
clearly extends to ${\bar V}^{\ell}$, and then restricts to the 
closure of ${\mathbb
Z}^d$ in ${\bar V}^{\ell}$.  It is not evident to me whether the 
$\var_y$'s for $y
\in {\mathbb Z}^d$ separate the points of this closure.  But even if 
they did, it is
not clear to me that we could then use this to apply the results of 
the previous
section to show that there are sufficient fixed-points in 
$\p_{\ell}{\mathbb Z}^d$
for the action of ${\mathbb Z}^d$.  It is this supply of fixed-points 
which we need
later.  So we take a more direct tack.  We show that every {\em 
linear} geodesic ray
in $V$ can be approximated by an almost-geodesic ray in ${\mathbb Z}^d$.  The
following lemma is closely related to Kronecker's theorem \cite{BrD}, 
so we just
sketch the proof.

\setcounter{lemma}{1}
\begin{lemma}
\label{lem7.2}
Let $v \in V$ with $\|v\| = 1$.  Then there is an unbounded strictly increasing
sequence $\{s_n\}$ of positive real numbers such that for every $\e > 
0$ there is an
$N$ such that if $s_n > N$ then there is an $x \in {\mathbb Z}^d$ for which
$\|x-s_nv\| < \e$.
\end{lemma}

\begin{proof}
If there is an $r \in {\mathbb R}^+$ with $rv \in {\mathbb Z}^d$ then 
we simply take
$s_n = nr$.  Suppose now that no such $r$ exists.  Consider the image 
of ${\mathbb
R}v$ in $V/{\mathbb Z}^d$.  Its closure is a connected subgroup, and 
so is a torus.
The dimension of this torus must be $\ge 2$ for otherwise there would 
be an $r$ as
above.  But for any finite closed interval $I$ of ${\mathbb R}$ the 
image of $Iv$ is
compact, and so must stay away from $0$ except at $0$.  Thus for any neighborhood 
of
$0$ there must be a $t$ outside of $I$ such that the image of $tv$ is in that
neighborhood.
\end{proof}

Let $\{s_n\}$ be as in the lemma.  Then we can find a subsequence, 
$\{t_k\}$, of the
sequence $\{s_n\}$, and for each $k$ we can choose a $x_k \in 
{\mathbb Z}^d$, such
that $\|x_k - t_kv\| < 1/k$ for all $k$.

\setcounter{lemma}{2}
\begin{lemma}
\label{lem7.3}
For $v$, $\{t_k\}$ and $\{x_k\}$ as above, define $\g$ by $\g(0) = 0$ 
and $\g(t_k) =
x_k$.  Then $\g$ is an almost-geodesic ray in $V$ which determines 
the same Busemann
point in $\p_{\ell}V$ as does the ray $t \mapsto tv$.
\end{lemma}

\begin{proof}
Given $\e > 0$, choose $N$ such that $1/N < \e/3$.  Then for $t_n \ge 
t_m \ge N$ we have from the triangle inequality
\begin{eqnarray*}
|\|x_n - x_m\| &+ &\|x_m\| - t_n| 
=  |\|x_n-x_m\| - \|(t_n-t_m)v\| + \|x_m\| - t_m| \\
&\le &\|(x_n-t_nv) - (x_m-t_mv)\| + \|x_m-t_mv\| < \e.
\end{eqnarray*}
 From this it follows that $\g$ is an almost-geodesic ray.  The fact that it
determines the same Busemann point as does $v$ now follows from Proposition
\ref{prop4.9}.
\end{proof}

\setcounter{proposition}{3}
\begin{proposition}
\label{prop7.4}
Let $v$ be a smooth point of the unit sphere of $V$, with support functional
$\s_v$.  Then there is a Busemann point $b_v \in \p_{\ell}{\mathbb 
Z}^d$ such that
for any $y \in {\mathbb Z}^d$ we have
\[
\var_y(b_v) = \s_v(y).
\]
Furthermore, $b_v$ is a fixed-point for the action of ${\mathbb Z}^d$ on
$\p_{\ell}{\mathbb Z}^d$.
\end{proposition}

\begin{proof}
Let $\g$ be an almost-geodesic ray associated with $v$ as in the 
above lemmas.  By
Proposition \ref{prop6.2} we know that
\[
\lim \var_y(x_k) = \s_v(y)
\]
for all $y \in V$.  But $\g$ is equally well an almost-geodesic ray 
in ${\mathbb
Z}^d$, and so defines a Busemann point $b_{\g} \in \p {\mathbb Z}^d$. 
But for $y
\in {\mathbb Z}^d$ its $\var_y$ for ${\mathbb Z}^d$ is just the restriction to
${\mathbb Z}^d$ of its $\var_y$ for $V$.  Thus $\var_y(b_{\g}) = 
\s_v(y)$ for $y \in
{\mathbb Z}^d$.  The proof that $b_{\g}$ is a fixed-point for the 
action is the same
as that for Proposition \ref{prop6.3}.
\end{proof}

We remark that, just as for $V$, different smooth points $v$ may have the same
$\s_v$, and so determine the same Busemann point of $\p_{\ell} 
{\mathbb Z}^d$, and so it can happen that
only a finite number of points of $\p_{\ell} {\mathbb Z}^d$ arise 
from smooth
points $v$.

We are now ready to prove one part of our Main Theorem 0.1, namely:

\setcounter{theorem}{4}
\begin{theorem}
\label{th7.5}
Let $\ell$ on ${\mathbb Z}^d$ be defined by $\ell(x) = \|x\|$ for a 
norm $\|\cdot\|$
on ${\mathbb R}^d$.  Let $L_{\ell}$ be defined on $C_c({\mathbb 
Z}^d,c)$ as before by
\[
L_{\ell}(f) = \|[M_{\ell},\pi_f]\|.
\]
Then $L_{\ell}$ is a Lip-norm on $C^*({\mathbb Z}^d,c)$.
\end{theorem}

\begin{proof}
Let $v$ be a smooth point of the unit sphere of $V$ for $\|\cdot\|$. 
Let $\s_v$
denote its support functional, and $b_v$ its corresponding Busemann 
point as above
in $\p_{\ell}{\mathbb Z}^d$.  Since $b_v$ is a fixed-point, it determines a
homomorphism from the cosphere algebra $C^*(G,C(\p_{\ell}G),\a,c)$ 
onto $C^*(G,c)$
which takes $M_{\var_y}$ to the constant $\s_v(y)$.  (We use here the 
amenability of
${\mathbb Z}^d$.)  Then under this homomorphism $[M_{\ell},\pi_f]$ is 
sent to the
operator
\[
\Sigma f(y)\var_y(b_v)\pi_y = \Sigma f(y)\s_v(y)\pi_y
\]
in $C^*({\mathbb Z}^d,c)$.  Let us denote this operator, and the
corresponding function, by $X_vf$. Of course
$\|X_vf\| \le L_{\ell}(f)$.

We let $\b$ denote the usual dual action \cite{Pdr} of the dual group 
${\hat G}$ on
$C^*(G,c)$ determined by
\[
(\b_p(f))(x) = \<x,p\>f(x)
\]
for $f \in C_c({\mathbb Z}^d)$ and $p \in {\hat G}$, where 
$\<\cdot,\cdot\>$ denotes
the pairing of $G$ and ${\hat G}$.  
Each $\tau \in V'$ determines an element of
${\hat G}$ by $\<x,\tau\> = \exp(i\tau(x))$ for $x \in {\mathbb 
Z}^d$.  Let $\G$
denote the lattice in $V'$ consisting of elements which on
${\mathbb Z}^d$ take values in $2\pi {\mathbb Z}$.  
Then we can identify ${\hat G}$ with the torus 
$V'/\G$, and then $V'$
is identified with the Lie algebra of ${\hat G}$, so that the 
exponential mapping is
just the quotient map from $V'$ to $V'/\G$.  The action $\b$ has an 
infinitesimal
version which is a Lie algebra homomorphism from the (Abelian) Lie 
algebra $V'$ into
the Lie algebra of derivations on $C^*(G,c)$.  
We denote it by $d\b$, and it is
determined by
\[
(d\b_{\tau}(f))(x) = i\tau(x)f(x).
\]
Each $f \in C_c(G)$ then determines a linear mapping, $\tau \mapsto 
d\b_{\tau}(f)$,
from $V'$ into $C^*(G,c)$, which we denote by $df$, much as done for 
theorem $3.1$
of \cite{R4}.

In terms of the notation just introduced, we see that for any smooth 
point $v$ we
have
\[
iX_vf = d\b_{\s_v}f = df(\s_v).
\]
With this notation our earlier inequality becomes
\[
\|df(\s_v)\| \le L_{\ell}(f).
\]
Now $V'$ has the dual norm $\|\cdot\|'$, and $C^*(G,c)$ has its 
$C^*$-norm.  So the
norm of the linear map $df$ between them is well-defined.  We denote it by
$\|df\|$.  But by Proposition \ref{prop6.7} the closed convex hull of 
the set of
$\s_v$'s is the unit ball in $V'$.  It follows that
\[
\|df\| \le L_{\ell}(f).
\]
But in theorem $3.1$ of \cite{R4} it is shown that $f \mapsto \|df\|$ is a
Lip-norm.  Thus we can apply comparison lemma $1.10$ of \cite{R4} to 
conclude that $L_{\ell}$ is a Lip-norm as well.
\end{proof}

%%%%%%%%%%%%%%%%%%%%%%%%%%%%%%%%%%%%%%%%%%%%%%%%%
\setcounter{section}{7}
\section{The boundary of $({\mathbb Z}^d,S)$}
\label{sec8}
%%%%%%%%%%%%%%%%%%%%%%%%%%%%%%%%%%%%%%%%%%%%%%%%%%

Let $S$ be a finite generating subset of $G = {\mathbb Z}^d$ such 
that $S=-S$ and $0
\notin S$.  Let $\ell$ denote the corresponding word-length function 
on $G$.  I do
not know how to give a concrete description of $\p_{\ell}G$.  (But note that
$\p_{\ell}G$ is totally disconnected since each $\var_y$ takes only 
integer values,
in contrast to what happens if $\ell$ comes, for example, from the 
Euclidean norm on
${\mathbb R}^d$.) We will show here how to construct a substantial supply of 
geodesic rays.
(Somewhat related considerations appear in \cite{Shp}, but geodesic rays and
compactifications are not considered there.)  In the next section we 
will show  that
our supply is sufficient to prove that when $M_{\ell}$ is used as the 
Dirac operator
for $C^*(G,c)$, then the corresponding metric on the state space of 
$C^*(G,c)$ gives
the weak-$*$ topology.

Our construction is motivated by several features which we found in Sections
\ref{sec6} and \ref{sec7}.  For convenience we view $G = {\mathbb 
Z}^d$ as embedded
in ${\mathbb R}^d$.
We let $K = K_S$ denote the (closed) convex hull in ${\mathbb R}^d$ of 
$S$.  Because $K$
is balanced (since $S = -S$), it determines a norm, $\|\cdot\|_S$, on ${\mathbb
R}^d$, for which it is the unit ball.  (In fact, $({\mathbb 
R}^d,\|\cdot\|_S)$ is
the ``asymptotic cone'' of $({\mathbb Z}^d,\ell)$---see exercise $8.2.12$ of 
\cite{BBI}.)  We
will see later that this norm is relevant.  The set of extreme points 
of $K_S$ is a
subset of $S$, which we will denote by $S^e$.  The faces of $K_S$ (of all
dimensions) will have certain subsets of $S^e$ as their extreme 
points, and will
intersect $S$ in certain subsets $F$.  Such an $F$ is characterized 
by the fact that
there is a linear functional $\s$ on ${\mathbb R}^n$ (not necessarily 
unique) such
that $\s(s) \le 1$ for all $s \in S$ and $F = \{s \in S: \s(s) = 
1\}$.  We call any
such $\s$ a {\em support functional} for $F$.  Note that $|\s(s)| \le 
1$ for all $s
\in S$.  By abuse of terminology we will refer to $F$ itself as a 
face of $K_S$, and
we will not distinguish between $F$ and the usual face which $F$ determines.

\setcounter{lemma}{0}
\begin{lemma}
\label{lem8.1}
Let $\s$ be a support functional for a face $F$ of $K_S$.  Then
\[
|\s(x)| \le \ell(x)
\]
for all $x \in G$.
\end{lemma}

\begin{proof}
Suppose that $x = \Sigma q(s)s$ for some function $q$ from $S$ to 
${\mathbb Z}$.
Then
\[
\s(x) = \Sigma q(s)\s(s) \le \Sigma |q(s)|.
\]
On considering the minimum for all such $q$, we see that $\s(x) \le 
\ell(x)$.  But
this holds for $-x$ as well, which gives the desired result.
\end{proof}

Let $F$ be a face of $K_S$.  Any function $\g$ from ${\mathbb Z}^+$ 
to $G$ which
consists of successively adding elements of $F$ (i.e., $\g(n+1) - 
\g(n) \in F$ for
$n \ge 0$) is a geodesic ray.  In fact, 
for any support functional $\s$ for $F$ the
above lemma tells us that we have $n \ge \ell(\g(n)) \ge \s(\g(n)) = 
n$.  Since $F$
is finite, some (perhaps all) elements of $F$ will have to be added 
in an infinite
number of times.  One can see that if the order in which the elements 
of $F$ are
added-in is changed, but the number of times they ultimately appear 
is the same,
then one obtains an equivalent geodesic ray.  A class of such geodesic 
rays can be specified
by a function on $F$ which has values either in ${\mathbb Z}^+$ or 
$+\infty$.  But
it seems to be tricky to decide when two such functions (possibly for 
different faces)
determine the same Busemann point.  For our present purposes we do not 
need to concern
ourselves with this issue.  It is sufficient for us to associate a canonical
geodesic ray to each face.  
This will be a special case of forming geodesic rays by
successively adding elements of the semigroup generated by $F$ (so 
that the domain
of the ray may be a proper subset of ${\mathbb Z}$).

\setcounter{notation}{1}
\begin{notation}
\label{nota8.2}
For a face $F$ of $K_S$ set $z_f = \Sigma\{s: s \in F\}$, and let 
$\g_F$ denote the
geodesic ray whose domain is $|F|{\mathbb Z}^+$ (where $|F|$ denotes 
the number of
elements of $F$) and which is defined by $\g(|F|n) = nz_F$.  We 
denote by $b_F$ the
corresponding Busemann point.  We denote by $G_F$ the subgroup of $G$ 
generated by
$F$.
\end{notation}

Again Lemma \ref{lem8.1} quickly shows that the above ray is geodesic.
The following proposition is analogous to Proposition 6.2.

\setcounter{proposition}{2}
\begin{proposition}
\label{prop8.3}
Let $\s$ be a support functional for a face $F$ of $K_S$.  For every 
$u \in G_F$ we
have
\[
\var_u(b_F) = \s(u).
\]
\end{proposition}

\begin{proof}
Since $u \in G_F$, there is a positive integer $N$ such that whenever 
$n \ge N$ then
$nz_F - u$ can be expressed as a sum of elements of $F$, so that 
$\ell(nz_F - u) =
\s(nz_F-u)$.  Of course $\ell(nz_F) = \s(nz_F)$.  Thus for $n \ge N$
\[
\var_u(nz_F) = \s(nz_F) - \s(nz_F-u) = \s(u).
\]
\end{proof}

\setcounter{proposition}{3}
\begin{proposition}
\label{prop8.4}
Let $F$ and $\s$ be as above.  For any $y \in G$ and $u \in G_F$ we have
\[
\var_{y+u}(b_F) = \var_y(b_F) + \s(u).
\]
\end{proposition}

\begin{proof}
Consider the set of $u$'s such that this equation holds for all $y 
\in G$.  It is easy to verify that 
this set is a subsemigroup of $G$.  But for $u$ in this set we have
\[
\var_{y-u}(\b_F) = \var_{(y-u)+u}(\b_F) - \s(u) = \var_y(\b_F) + \s(-u),
\]
so that this set is a group.  It thus suffices to verify the above 
equation for each
$u = s \in F$.

So let $s \in S$.  Since $n \mapsto \var_y(nz_F)$ is integer-valued, 
non-decreasing by Lemma 4.5, and
bounded, we can find a positive integer $N$ such that
\[
\var_y(b_F) = \ell((N+m)z_F) - \ell((N+m)z_F-y)
\]
for all $m \ge 0$.  We can find a larger $N$ such that also
\[
\var_{y+s}(b_F) = \ell((N+m)z_F) - \ell((N+m)z_F - (y+s))
\]
for all $m \ge 0$.  Since $\s(s) = 1$ it is then clear that we need 
to show that
\[
\ell((N+m)z_F - (y+s)) = \ell((N+m)z_F-y) - 1
\]
for some $m \ge 0$.  Let ${\bar y} = y - Nz_F$.  Then what we need becomes
\[
\ell(mz_F - ({\bar y} + s)) = \ell(mz_F - {\bar y}) - 1
\]
for some $m \ge 0$.  Note that $\ell(mz_F) - \ell(mz_F-{\bar y})$ is 
independent of
$m \ge 0$ because of our choice of $N$, and similarly for ${\bar y} + 
s$ instead of
${\bar y}$.

Since $S = -S$ and $0 \notin S$, we can find a subset, $S^+$, such 
that $S^+ \cup (-S^+) = S$ and $S^+ \cap (-S^+) = \emptyset$.  Since $F \cap (-F) = 
\emptyset$, we can
require that $F \subseteq S^+$.  Index the elements of $S^+$ in such 
a way that $s_1
= s$, and $F = \{s_1,\dots,s_{|F|}\}$, where $|F|$ denote the number 
of elements in
$F$.  Since $S$ generates $G$, we can express ${\bar y}$ as ${\bar y} = \Sigma
n_js_j$ where $n_j \in {\mathbb Z}$ for each $j$.  Then $\ell({\bar 
y})$ will be the
minimum of the sums $\Sigma |n_j|$ over all such expressions for 
${\bar y}$.  We
make a specific choice of such a minimizing set $\{n_j\}$.  (It need 
not be unique.)

Since $\ell(mz_F) = m|F|$ by Lemma \ref{lem8.1}, the stability 
described earlier
says that $m|F| - \ell(mz_F-{\bar y})$ is independent of $m \ge 0$. 
We combine this
for $m=0$ and $m=1$ to obtain $-\ell(-{\bar y}) = |F| - 
\ell(z_F-{\bar y})$.  We use
this to calculate
\begin{eqnarray*}
|F| &+ &\Sigma|n_j| = |F| + \ell(-{\bar y}) = \ell(z_F - {\bar y}) \\
&= &\ell\left( \sum_{j \le |F|} (1-n_j)s_j + \sum_{j > |F|} n_js_j\right) 
\le \sum_{j \le |F|} |1-n_j| + \sum_{j>|F|} |n_j|.
\end{eqnarray*}
On comparing the two ends, we see that we must have $n_j \le 0$ for 
$j \le |F|$, and
that the two ends must be equal.  Thus
\[
\ell(z_F - {\bar y}) = \sum_{j \le |F|} (1-n_j) + \sum_{j > |F|} |n_j|.
\]
Now
\begin{eqnarray*}
z_F - {\bar y} - s &= &\sum_{j \le |F|} (1-n_j)s_j + \sum_{j>|F|} 
n_js_j - s_1 \\
&= &-n_1s_1 + \sum_2^{|F|} (1-n_j)s_j + \sum_{j>|F|} n_js_j.
\end{eqnarray*}
 From the fact that $n_j \le 0$ for $j \le |F|$ it follows that
\begin{eqnarray*}
\ell(z_F - {\bar y} - s) &\le &-n_1 + \sum_2^{|F|} (1-n_j) + 
\sum_{j>|F|} |n_j| \\
&= &-1 + \ell(z_F - {\bar y}).
\end{eqnarray*}
 From the triangle inequality and the fact that $\ell(s) = 1$ it follows that
\[
\ell(z_F - {\bar y} - s) = -1 + \ell(z_F - {\bar y}),
\]
as needed.
\end{proof}

\setcounter{corollary}{4}
\begin{corollary}
\label{cor8.5}
For any $y,z \in G$ and any $u \in G_F$, and for any support 
functional $\s$ for
$F$, we have
\[
\var_{y+u}(\a_z(b_F)) = \var_y(\a_z(b_F)) + \s(u).
\]
\end{corollary}

\begin{proof}
Using the $1$-cocycle identity $2.2$ and Proposition \ref{prop8.4} we obtain
\begin{eqnarray*}
\var_{y+u}(\a_z(b_F)) &= &(\a_{-z}\var_{y+u})(b_F) 
= \var_{y-z+u}(b_F) - \var_{-z}(b_F) \\
&= &\var_{y-z}(b_F) + \s(u) - \var_{-z}(b_F) \\
&= &(\a_{-z}\var_y)(b_F) + \s(u) = \var_y(\a_z(b_F)) + \s(u).
\end{eqnarray*}
\end{proof}

\setcounter{proposition}{5}
\begin{proposition}
\label{prop8.6}
Let $F$ be a face of $K$.  For each $u \in G_F$ the homeomorphism 
$\a_u$ of ${\bar
G}^{\ell}$ leaves fixed each point of the $\a$-orbit of $b_F$.  That 
is, for each $z
\in G$ we have
\[
\a_u(\a_z(b_F)) = \a_z(b_F).
\]
\end{proposition}

\begin{proof}
Because $G$ is Abelian, it suffices to show that $\a_u(b_F) = b_F$. 
For this we
must verify that $f(\a_u(b_F)) = f(b_F)$ for all $f \in C({\bar 
G}^{\ell})$.  It
suffices to verify this for $f = \var_y$ for each $y \in G$.  But from the
$1$-cocycle identity 2.2 and Proposition \ref{prop8.4} we have
\begin{eqnarray*}
\var_y(\a_u(b_F)) &= &(\a_{-u}\var_y)(b_F) = \var_{y-u}(b_F) - 
\var_{-u}(b_F) \\
&= &\var_y(b_F) + \s(-u) + \s(u) = \var_y(b_F).
\end{eqnarray*}
\end{proof}

We will also need the following fact:

\setcounter{proposition}{6}
\begin{proposition}
\label{prop8.7}
If $y \notin G_F$ then $\var_y$ is not constant on the $G$-orbit of 
$b_F$, and in
fact there is an $s \in S$ such that $s \notin F$ and
\[
\var_y(\a_s(b_F)) = \var_y(b_F) + (1-\var_{-s}(b_F)),
\]
with $\var_{-s}(b_F) = 0$ or $-1$.
\end{proposition}

\begin{proof}
Let $S^+$ and the indexing $\{s_j\}$ be as in the proof of Proposition
\ref{prop8.4}.  Much as in that proof, we can find a large enough $N$ 
that $\var_{y
\pm s_j}((N+m)z_F)$ is constant for $m \ge 0$ for all $\pm s_j$ 
simultaneously, as
is $\var_y((N+m)z_F)$.  Set ${\bar y} = y - Nz_F$.  For this ${\bar y}$ choose
$\{n_j\}$ as before so that ${\bar y} = \Sigma n_js_j$ and 
$\ell({\bar y}) = \Sigma
|n_j|$.  Since $y \notin G_F$, also ${\bar y} \notin G_F$, and so 
there is a $k >
|F|$ such that $n_k \ne 0$.  Suppose that $n_k \ge 1$.  Then
\[
{\bar y} - s_k = \sum_{j \ne k} n_js_j + (n_k-1)s_k,
\]
so that
\[
\ell({\bar y} - s_k) \le \sum_{j \ne k} |n_j| + n_k-1 = \ell({\bar y}) - 1.
\]
 From the triangle inequality we then obtain $\ell({\bar y} - s_k) = 
\ell({\bar y}) - 1$, that is,
\[
\ell(Nz_F - y + s_k) = \ell(Nz_F - y) - 1.
\]
 From our choice of $N$ (and with $m=0$) we then get
\begin{eqnarray*}
\var_{y-s_k}(b_F) &= &\var_{y-s_k}(Nz_F) \\
&= &N|F| - \ell(Nz_F - y + s_k) = N|F| -
\ell(N|F| - y) + 1 \\
&= &\var_y(b_F)+1.
\end{eqnarray*}
We combine this with the $1$-cocycle identity $2.2$ to obtain
\begin{eqnarray*}
\var_y(\a_{s_k}(b_F)) &= &(\a_{-s_k}\var_y)(b_F) \\
&= &\var_{y-s_k}(b_F) - \var_{-s_k}(b_F) \\
&= &\var_y(b_F) + (1-\var_{-s_k}(b_F)).
\end{eqnarray*}

Since $\var_{-s_k}$ takes only the values $0$, $\pm 1$, the desired 
conclusion is
then obtained from:

\setcounter{lemma}{7}
\begin{lemma}
\label{lem8.8}
If $s \in S$ and $\var_s(b_F) = 1$ then $s \in F$.
\end{lemma}

\begin{proof}
If $\var_s(b_F) = 1$, then for large $n$, and for a support 
functional $\s$ for $F$,
we have
\begin{eqnarray*}
n|F| - 1 &= &\ell(nz_F) - 1 = \ell(nz_F - s) \\
&\ge &\s(nz_F-s) = n|F| - \s(s),
\end{eqnarray*}
so that $1 \le \s(s)$, and so $s \in F$.
\end{proof}

The above argument for the proof of Proposition \ref{prop8.7}
was under the assumption that $n_k \ge 1$.  If 
instead we have
$n_k \le -1$, then we carry out a similar argument using $-s_k$ 
instead of $s_k$.
This concludes the proof of Proposition \ref{prop8.7}.
\end{proof}

%%%%%%%%%%%%%%%%%%%%%%%%%%%%%%%%%%%%%%%%%%%%%%
\setcounter{section}{8}
\section{Word-length functions give Lip-norms on $C^*({\mathbb Z}^d,c)$}
\label{sec9}
%%%%%%%%%%%%%%%%%%%%%%%%%%%%%%%%%%%%%%%%%%%%%%%%%%%%

We will now see how to use the results of the previous section to 
prove the part of
our Main Theorem \ref{th0.1} concerning word-length functions.  We 
use the notation
of the previous section, and in particular, the norm $\|\cdot\|_S$ 
determined by $K
= K_S$.  Here we will consider the (proper) faces of $K$ of maximal dimension,
namely of dimension $d-1$.  We will call them ``facets'' of $K$, as is not
infrequently done.  The interior points of the facets are the smooth 
points of the
unit sphere for $\|\cdot\|_S$.  Again our terminology and notation will not
distinguish between facets as intersections 
of $K$ with hyperplanes, and as the
corresponding subsets of $S$.  Because $K$ has only a finite number of extreme
points, every point of the boundary of $K$ is contained in at least 
one facet, and
there are only a finite number of facets.  Each facet $F$ has a unique support
functional, which we denote by $\s_F$.  Furthermore, $F$ contains a basis for
${\mathbb R}^d$, and consequently $G_F$ is of finite index in $G$. 
This has the
crucial consequence for us that the orbit, ${\mathcal O}_F$, of $b_F$ in
$\p_{\ell}G$ under the action $\a$, is finite.  (Apply Proposition 
\ref{prop8.6}.)
We consider the restriction map from $C(\p_{\ell}G)$ onto $C({\mathcal O}_F)$.
Since it is $\a$-equivariant, 
it gives an algebra homomorphism, $\Pi_F$, from
$C^*(G,C(\p_{\ell}G),\a,c)$ onto $C^*(G,C({\mathcal O}_F),\a,c)$.  If 
we let $\pi$
and $M$ denote also the corresponding homomorphisms of $G$ and 
$C({\mathcal O}_F)$
into this latter algebra, and if for each $y \in G$ we let $\psi_y$ denote the
restriction of $\var_y$ to ${\mathcal O}_F$, then
\[
\Pi_F([M_{\ell},\pi_f]) = \Sigma f(y)M_{\psi_y}\pi_y.
\]
Let $Q$ be a set of coset representatives for $G_F$ in $G$ containing $0$. 
Then we can express the
above as
\[
\Sigma_{q \in Q} (\Sigma_{u \in G_F} f(u+q)M_{\psi_{u+q}}{\bar 
c}(u,q)\pi_u)\pi_q.
\]
From Corollary \ref{cor8.5} we see that $\psi_{u+q} = \psi_q + 
\s_F(u)$.  For each
$q$ let $g^q$ be the function on $G_F$ defined by $g^q(u) = 
f(u+q){\bar c}(u,q)$.
We can also view $g^q$ as a function on $G$ by giving it value $0$ 
off $G_F$.  Then
we can rewrite our previous expression for $\Pi_F([M_{\ell},\pi_f])$ as
\[
\Sigma_q(\Sigma_u g^q(u)(\s_F(u) + M_{\psi_q})\pi_u)\pi_q.
\]
As in Section \ref{sec7} let ${\hat G} = {\mathbb T}^d$ be the dual 
group of $G$,
and denote the pairing between $G$ and ${\hat G}$ by $\<x,s\>$.  Let 
$\b$ now denote
the usual dual action of ${\hat G}$ on $C^*(G,C({\mathcal O}_F),\a,c)$, so that
\[
\b_s(M_{\psi}\pi_x) = \<x,s\>M_{\psi}\pi_x.
\]
Then the finite group $(G/G_F)^{\wedge}$ can be identified with the set of
characters on $G$ which take value $1$ on $G_F$.  We can thus restrict $\b$ to
$(G/G_F)^{\wedge}$ and average over $(G/G_F)^{\wedge}$.   
This gives a projection of norm~$1$ onto the subalgebra 
of elements
supported on $G_F$, and this projection on functions on $G$ is just 
restriction of
functions to $G_F$.  If for each fixed $q$ we apply this projection to 
the product with $\pi_q^*$
of the above expression for $\Pi_F([M_{\ell},\pi_f])$, we find that
\[
\|[M_{\ell},\pi_f]\| \ge \|\Sigma_u g^q(u)(\s_F(u)+M_{\psi_q})\pi_u\|
\]
for each $q$.  The norm on the right is that of 
$C^*(G,C({\mathcal O}_F),\a,c)$.  But section 2.27 of \cite{ZM}
tells us that $C^*(G_F,C({\mathcal O}_F),\a,c)$ is a $C^*$-subalgebra
of $C^*(G,C({\mathcal O}_F),\a,c)$ under the evident identification
of functions.
Thus we can view the operator on the right 
as being in $C^*(G_F,C({\mathcal O}_F),\a,c)$, where we are here
restricting $\a$ and $c$ to $G_F$.  But from Proposition 
\ref{prop8.6} we know that
the action $\a$ of $G_F$ on ${\mathcal O}_F$ is trivial.  Thus we have the
decomposition
\[
C^*(G_F,C({\mathcal O}_F),\a,c) \cong C({\mathcal O}_F) \otimes C^*(G_F,c).
\]
Let $a_q = \Sigma g^q(u)\pi_u$ and $b_q = \Sigma g^q(u)\s_F(u)\pi_u$. 
Then in terms
of the above decomposition we are looking at $I \otimes b_q + \psi_q 
\otimes a_q$.
 From Proposition \ref{prop8.7} we know that $\psi_q$ is not constant 
on ${\mathcal
O}_F$ for $q \ne 0$.  Note that $\psi_0 \equiv 0$.  For given $q \ne 
0$ let $m_j$
for $j = 1,2$ be two distinct values of $\psi_q$.  Upon evaluating at 
the points
where $\psi_q$ takes these values, and using our earlier inequality, 
we see that
\[
\|b_q + m_ja_q\| \le \|[M_{\ell},\pi_f]\| = L_{\ell}(f)
\]
for $j = 1,2$.  Upon writing the inequalities as
\[
\|m_j^{-1}b_q + a_q\| \le |m_j|^{-1}L_{\ell}(f)
\]
and using the triangle inequality to eliminate $a_q$, and 
simplifying, we find that
\[
\|b_q\| \le (|m_1| + |m_2|)/|m_1-m_2|L_{\ell}(f).
\]
(If either $m_j$ is $0$ the path is simpler.)  Of course $m_1$ and 
$m_2$ depend on
$q$.  Thus we see that we have found a constant, $k_q$, 
such that $\|b_q\| \le
k_qL_{\ell}(f)$.  For $q = 0$ we have the same inequality with $k_0 = 
1$ since $\psi_0 = 0$.  Much as in Section \ref{sec7} set $X_Ff = \Sigma 
\s_F(x)f(x)\pi_x$.  Then
\begin{eqnarray*}
X_Ff &= &\Sigma \s_F(x)f(x)\pi_x 
= \Sigma_q(\Sigma_{u \in G_F} \s_F(u+q)f(u+q){\bar c}(u,q)\pi_u)\pi_q \\
&= &\Sigma_q \s_F(q) (\Sigma_u\s_F(u)g^q(u)\pi_u)\pi_q 
= \Sigma \s_F(q)b_q\pi_q.
\end{eqnarray*}
When we combine this with the inequality obtained earlier for 
$\|b_q\|$, we obtain
\[
\|X_Ff\| \le (\Sigma |\s_F(q)|k_q)L_{\ell}(f).
\]
Observe that the $\s_F(q)$'s and $k_q$'s do not depend on $f$, but 
only on $F$ and
the choice $Q$ of coset representatives.  Thus for each facet $F$ we 
have obtained a
constant, $k_F$, such that
\[
\|X_Ff\| \le k_FL_{\ell}(f)
\]
for all $f \in C_c(G)$.  Note that knowing that $k_F$ is finite
is the crucial place where we use that the number of coset
representatives in $Q$ is finite.

Just as toward the end of Section \ref{sec7}, we have the dual action $\b$ of
${\mathbb T}^d$ on $C^*(G,c)$, and the corresponding differential 
$df$ of any $f \in
C_c(G)$, such that $df(\s_F) = iX_Ff$.  Then our inequality above 
gives, much as in
Section \ref{sec7},
\[
\|df(\s_F)\| \le k_FL_{\ell}(f).
\]

Recall now the norm $\|\cdot\|_S$ determined by $K = K_S$.  The $\s_F$'s are
exactly the support functionals corresponding to the smooth points of the unit
sphere for $\|\cdot\|_S$.  Let $\|df\|_S$ denote the norm of the 
linear map $df$
using the dual norm $\|\cdot\|'_S$.  Also let $k = \max\{k_F: F 
\mbox{ is a facet}\}$.
Then from Proposition \ref{prop6.7} we conclude, much as in Section 
\ref{sec7}, that
\[
\|df\|_S \le kL_{\ell}(f).
\]
Then just as in Section \ref{sec7} we conclude that $L_{\ell}$ is a Lip-norm.
This concludes the proof of Main Theorem 0.1. \\  \qed

Since the norm $\|\cdot\|'_S$ on $V'$ does not come from an inner product, 
and $V'$ can be thought of as the analogue of the tangent space
at the non-existent points of the quantum space $C^*(G, c)$, we can consider
that we have here a non-commutative Finsler geometry (as also in section 3
of \cite{R4}). The metric geometry from $L_\ell$ also, in a vague
way, seems Finsler-like.

I imagine that the above considerations can be generalized so that the 
Main Theorem
can be extended to weighted-word-length functions, where each 
generator has been
assigned a weight.  I imagine that they can also be generalized
to deal with extensions of ${\mathbb Z}^d$ by finite groups.
But I have not explored these possibilities.

Since our estimates for the proof of the Main Theorem depend just 
on the behavior of the $\var_y$'s on the boundary, the conclusions of the 
Main Theorem will also be valid if $\ell$ is replaced by the
translation-bounded function $\ell + h$ where $h$ is any function 
in $C_\infty({\mathbb Z}^d)$. 

\setcounter{section}{9}
%%%%%%%%%%%%%%%%%%%%%%%%%%%%%%%%%%%%%%%
\section{The free group}
%%%%%%%%%%%%%%%%%%%%%%%%%%%%%%%%%%%%%%%%
\label{sec10}

We briefly discuss here how the ideas developed earlier apply to the free
(non-Abelian) group on two generators, $G = F_2$.  Denote the two 
generators by $a$
and $b$, and take them and their inverses as our generating set $S$. 
Let $\ell$
denote the corresponding length function.  It is well-known 
\cite{GhH} that $F_2$ is
a hyperbolic group, and that its Gromov boundary, $\p_hG$, is 
described as the set
of all infinite (to the right) reduced words in the elements of $S$. 
(The ``$h$''
in $\p_hG$ is for ``hyperbolic''---it does not denote a length function.)  The
action of $G$ on $\p_hG$ is the evident one 
by ``left concatenation'' (and then
reduction).  We can obtain the topology of $\p_hG$ and of the 
compactification of
$G$ as follows.  (See comment ii) on page~$104$ of \cite{GhH}.)  
To include the
elements of $G$ we need a ``stop'' symbol.  We denote it by $p$.  We 
let $S'$ denote
$S$ with $p$ added, and we let $\stackrel{\infty}{\prod} S'$ denote the set of
sequences with values in $S'$, with its compact topology of ``index-wise''
convergence.

\setcounter{notation}{0}
\begin{notation}
\label{nota10.1}
Let ${\bar G}^h$ be the subset of $\stackrel{\!\!\infty}{\prod} S'$ 
consisting of all
sequences such that
\begin{itemize}
\item[1)] If $p$ occurs in the sequence then all subsequent letters in that
sequence are $p$.
\item[2)] The sequence is reduced, in the sense that $a$ and $a^{-1}$ are
never adjacent entries, and similarly for $b$ and $b^{-1}$.
\end{itemize}
\end{notation}

It is easily seen that ${\bar G}^h$ is a closed subset of
$\stackrel{\infty}{\prod} S'$, so compact.  We identify the elements of 
$G$ with the
words containing $p$ (and in particular, we identify the identity 
element of $G$
with the constant sequence with value $p$).  With this understanding, 
it is easily
seen that $G$ is an open dense subset of ${\bar G}^h$.  We identify 
$\p_hG$ with the
infinite words which do not contain $p$.

The group $G$ again acts on ${\bar G}^h$ by left concatenation.  It is 
easily seen that
this action is by homeomorphisms.  Consider the function $\var_a$ on 
$G$.  For any
word $w$ we have $\ell(a^{-1}w) = \ell(w) + 1$ if $w$ begins with the 
letters $a^{-1}$,
$b$ or $b^{-1}$, or is the identity element, while $\ell(a^{-1}w) = 
\ell(w) - 1$ if
$w$ begins with the letter $a$.  
Thus $\var_a(w) = \ell(w) - \ell(a^{-1}w)$ has
value $1$ if $w$ begins with the letter $a$, and value $-1$ 
otherwise.  But we can
extend $\var_a$ to ${\bar G}^h$ by exactly this same prescription, 
and it is easily
seen that this extended $\var_a$ is continuous on ${\bar G}^h$.  We 
do the same with
$\var_b$, $\var_{a^{-1}}$ and $\var_{b^{-1}}$.  By using the 
$1$-cocycle identity
$2.2$ inductively, we see that each $\var_x$ for $x \in G$ extends to 
a continuous
function on ${\bar G}^h$ (in a unique way since $G$ is dense).  Of course the
functions in $C_{\infty}(G)$ extend by giving them value $0$ on 
$\p_hG$, and the
constant functions also extend.  In this way we identify $C({\bar 
G}^{\ell})$ with a
unital subalgebra of $C({\bar G}^h)$.

Let us see now that the subalgebra $C({\bar G}^{\ell})$ 
separates the points of
${\bar G}^h$.  Because the subalgebra contains $C_{\infty}(G)$, it is 
clear that we
only need to treat the points of $\p_hG$.  Let $v,w \in \p_hG$ with 
$v \ne w$.  Then
there must be a first entry where they differ.  That is, we can write 
them as $v =
x{\tilde v}$, $w = x{\tilde w}$ where $x$ is a finite word while 
${\tilde v}$ and
${\tilde w}$ differ in their first entry.  Suppose the first entry of 
${\tilde v}$
is $a$ while the first entry of ${\tilde w}$ is not $a$.  Then from 
what we saw above
\[
(\a_x\var_a)(v) = \var_a(x^{-1}v) = \var_a({\tilde v}) = 1,
\]
while in the same way $(\a_x\var_a)({w}) = -1$.  Thus the 
subalgebra $C({\bar
G}^{\ell})$ separates the points of ${\bar G}^h$, and so by the 
Stone--Weierstrass
theorem $C({\bar G}^{\ell}) = C({\bar G}^h)$, so that ${\bar G}^{\ell} = {\bar
G}^h$.  Thus in this case the metric and hyperbolic boundaries coincide. (The
referee has pointed out that if instead we take as generating set 
$\{a^{\pm 1}, a^{\pm 2}, b^{\pm 1}\}$, then the resulting metric
compactification will be different from that described above, because
just as in Example 5.2 we will obtain two ``parallel" geodesic rays,
namely $(e, a^2, a^4, \dots)$ and $(a, a^3, a^5, \dots)$, which will
give different Busemann points.) 

Each $w \in \p_hG$ specifies a unique geodesic ray to it from $e$, namely
$e,w_1,w_1w_2,w_1w_2w_3,\dots$.  
Thus every point of $\p_{\ell}G$ is a Busemann
point.  It is well-known \cite{AD2} that the action of $G$ on $\p_hG$ 
is amenable.
If one uses the definition of amenability in terms of maps from $\p_hG$ to
probability measures on $G$ which was stated in Section \ref{sec3}, 
then this is
seen by letting the $n$-th map, $m_n$, be the map which assigns to $w 
\in \p_hG$ the
probability measure which gives mass $1/n$ to the first $n$ points of 
the geodesic
ray from $e$ to $w$ \cite{AD2}.  In view of Theorem \ref{th3.7} this 
implies that
the cosphere algebra $S_{\ell}^*A$ for the spectral triple $(A = 
C_r^*(F_2),\bell^2(F_2),M_{\ell})$ is
$C^*(G,C(\p_hF_2),\a)$.

However, the action $\a$ on $\p_hF_2$ does not have any finite 
orbits, and so I do
not see how to continue along the lines of the previous section to determine
whether the metric on the state space $S(C_r^*(F_2))$ coming from the above
spectral triple gives the state space the weak-$*$ topology, or even 
just finite
diameter.  The difficulty remains:  What information can one obtain about
$\|\pi_f\|$ if one knows that $\|[M_{\ell},\pi_f]\| \le 1$?

%\bibliographystyle{plain}

%\bibliography{references}

\end{document}